\documentclass{elsarticle}

\usepackage{setspace}
\usepackage{amsmath}
\usepackage{amsfonts}
\usepackage{mathrsfs}
\usepackage{bm}
\usepackage{bbm}
\usepackage{graphicx}
\usepackage{tikz}
\usepackage{pgfplots}
    \usetikzlibrary{spy}
\usepackage{caption}
\usepackage{subcaption}
\pgfplotsset{compat = 1.3}
\usepackage[margin=1in]{geometry}
\usepackage[ruled,vlined]{algorithm2e}
\usepackage{xcolor}

\journal{Journal of Computational Physics}

\begin{document}

\begin{frontmatter}

\title{A Face-Upwinded Spectral Element Method }

\author[rvt2,rvt3]{Y.~Pan\fnref{fn1}\corref{cor1}}
\ead{yllpan@berkeley.edu}

\author[rvt2,rvt3]{P.-O.~Persson\fnref{fn3}}
\ead{persson@berkeley.edu}

\address[rvt2]{Department of Mathematics, University of California, Berkeley, Berkeley, CA 94720, United States}
\address[rvt3]{Mathematics Group, Lawrence Berkeley National Laboratory, 1 Cyclotron Road, Berkeley, CA 94720, United States}
\cortext[cor1]{Corresponding author}
\fntext[fn1]{Graduate student, Department of Mathematics, University of California, Berkeley}
\fntext[fn3]{Professor, Department of Mathematics, University of California, Berkeley}


\begin{abstract}
We present a new high-order accurate discretisation on unstructured meshes of quadrilateral elements. Our Face Upwinded Spectral Element (FUSE) method uses the same node distribution as a high-order continuous Galerkin (CG) method, but with a particular choice of node locations within each element and an upwinded stencil on the face nodes. This results in a number of benefits, including fewer degrees of freedom and straight-forward integration with CG. We present the derivation of the scheme and the analysis of its properties, in particular showing stability using von Neumann analysis. We show numerical evidence for its accuracy and efficiency on multiple classes of problems including convection-dominated flows, Poisson's equation, and the incompressible Navier-Stokes equations.
\end{abstract}

\end{frontmatter}
	
\section{Introduction}

Over the last few decades, significant research has been dedicated to the development of stable, high-order accurate numerical methods for convection-dominated flow problems. While the continuous Galerkin finite element method is a popular approach, and often used together with collocation-based spectral element methods \cite{karniadakis05spectral,shen11spectral}, it requires specialised stabilisation techniques. Methods based on artificial diffusion such as SUPG, VMS, and spectrally vanishing viscosity \cite{tadmor89spectral,deville2002high,hughes2018multiscale} have been proposed but for various reasons can be difficult to generalise to very complex problems. Consequently, researchers have developed several methods based on discontinuous solution fields, with the discontinuous Galerkin (DG) method \cite{Reed_Hill,cockburn01rkdg,hesthaven08dgbook} being the most popular. It offers provable linear stability for any polynomial degree and element shape. Related methods include DG-SEM \cite{kopriva96,gassner10dgsem}, spectral differences (SD) and spectral volumes (SV) \cite{zj02spectralvolume,zj06spectraldifference,wang2016spectral}, and flux reconstruction (FR) \cite{huynh2007flux}, which can be shown to be identical in some special cases \cite{zj14comparison}, but in general define different schemes with varying numerical properties.

Here, we introduce the Face-Upwinded Spectral Element (FUSE) method, which aims to combine the ideas from the discontinuous methods to obtain stabilisation on continuous solution fields. The approach is straightforward: we employ a standard spectral differentiation technique for all interior nodes, while an upwinded high-order stencil is utilised only for the nodes on element faces. This methodology, along with a unique set of node locations, results in a linearly stable scheme for any degree. The primary motivation behind this method is in its simplicity, both due the fewer degrees of freedom compared to DG and in ease of assembly. However, we observe other advantages, such as superior CFL conditions and the potential for the use of improved solvers such as static condensation, which is less clear how to apply on discontinuous methods.

We first describe the method in detail for the one dimensional case. Using von Neumann analysis, we demonstrate that using traditional node choices such as equidistant or Gauss-Lobatto the scheme is actually unstable already for cubic approximations. However, stability can be achieved using an unusual node choice of Gauss-Legendre plus boundary points. We show that for constant-coefficient problems, the scheme can be re-written as a Spectral Difference method with a specific choice of solution and flux nodes. Furthermore we also show that in general it can be understood as a nodally integrated Petrov-Galerkin method with the specific choice of nodes. This directly shows convergence for arbitrary polynomial degrees using these nodes, by applying previous stability results. We also show how to discretise second order operators using an upwind-downwind strategy similar to the Local Discontinuous Galerkin (LDG) method \cite{ldg1998}. Next, we extend the scheme to higher space dimensions, where in particular we discuss how upwinding is performed on boundary nodes. Our numerical examples show evidence of the high-order accuracy for a range of problems, including convection, diffusion, and the incompressible Navier-Stokes equations.

The method is closely related to many previously proposed numerical schemes, and it can be argued that it is only a minor modification of several of the methods mentioned above. However, we are not aware of any other work using these particular choices, which are critical to obtain the attractive properties of our scheme. As discussed above, the method can be shown to be identical to the SD Method \cite{zj06spectraldifference} for a special case. However, in general it uses different solution nodes and continuous solutions which leads to very different properties. The unusual node choices were also used by Jameson \cite{jameson2009} as flux nodes in the SD method, but as far as we know there are no other schemes that are based on using these nodes as solutions nodes (continuous or discontinuous). The extensions to 2D also have many similarities with techniques used in the finite difference community \cite{leveque07finite}, but we note that due to the unstructured meshes the details end up being quite different. 

\section{1st derivative operators in 1D} \label{sect:conv1d}
\subsection{Preliminaries}
For this section we consider the general first order equation in conservative form
\begin{equation} \label{eq:conserv1d}
    \frac{\partial u}{\partial t} + \frac{\partial F(u)}{\partial x} = 0
\end{equation}
on the domain $\Omega = [0,1]$ with periodic boundary conditions. Assuming sufficient continuity on the flux function $F$ using the chain rule this can be rewritten as 
\begin{equation}
    \label{eq:non_conserv_convect}
    \frac{\partial u}{\partial t} + a(u) \cdot \frac{\partial u}{\partial x} = 0
\end{equation}
where $a(u) = F'(u)$. This form is useful in the formulation of upwind methods, where $a(u)$ is taken to be the velocity in constructing upwind discretisations.

\subsection{Discretisation}
\begin{figure}
    \centering
    \includegraphics[scale=0.4]{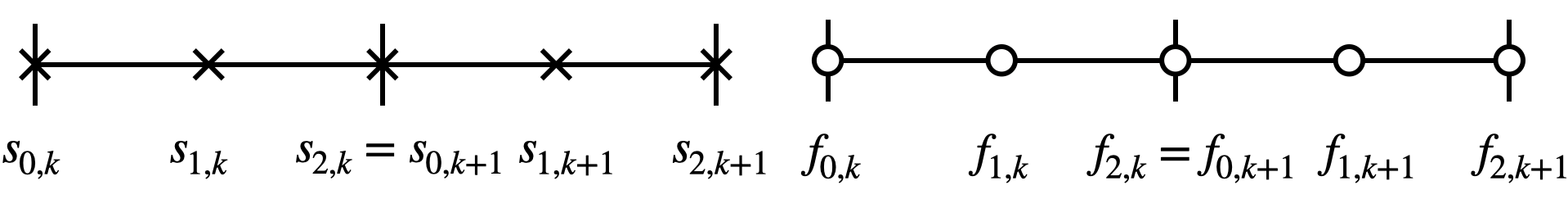}
    \caption{Example discretisation of 1D domain with two elements for $p=q=2$. Solution nodes are denoted with crosses and flux nodes with circles. Both solution and flux nodes at element boundaries are not repeated.}
    \label{fig:discretisation1d}
\end{figure}

To discretise the solution $u$ the domain $\Omega$ is paritioned into distinct elements $\mathcal{T}_h$ analogously to the Finite Element (FEM) or Discontintuous Galerkin (DG) Methods. Two sets of nodes are distributed within each element $K \in \mathcal{T}_h$, a set of solution nodes $\{ s_{0}, ..., s_{p} \}$ which are used to discretise the solution $u(x)$, and a set of flux nodes $\{ f_{0}, ..., f_{q} \}$, where $q \geq p$, for discretising the flux $F(u)$. To ensure coupling between elements the solutions nodes $s_{0}, s_{p}$ and flux nodes $f_{0},f_{q}$ are placed at the endpoints of the element. However unlike with DG and other methods involving a discontinuous solution field, solution and flux nodes on element boundaries are not repeated as shown in Fig. \ref{fig:discretisation1d}.

We introduce the function spaces $V_s(\mathcal{T}_h), V_f(\mathcal{T}_h)$ on elements $\mathcal{T}_h$ for the solution and flux fields respectively as follows
\begin{align}
    V_s(\mathcal{T}_h) &= \{ v \in H^1(\Omega) : v|_K \in \mathcal{P}_{p},~ \forall K \in \mathcal{T}_h  \} \\
    V_f(\mathcal{T}_h) &= \{ v \in H^1(\Omega) : v|_K \in \mathcal{P}_{q},~ \forall K \in \mathcal{T}_h \}
\end{align}
where $\mathcal{P}_{t}$ denotes the space of polynomials of degree at most $t \geq 1$. In practice a set of interpolating polynomials $\{ \phi_{i}^s \}, \{ \phi_{i}^f \}$ such that 
\begin{equation}
    \phi_{i}^s(s_{j}) = \delta_{ij}, ~\phi_{i}^f(f_{j}) = \delta_{ij}
\end{equation} 
are chosen as basis sets for the two spaces respectively. These are the standard interpolating basis functions used commonly in FEM. Continuity of functions in these spaces follow as a consequence of the shared boundary nodes across elements as with standard FEM. The basis function $\phi_{i}^s$ is said to be associated with solution node $s_{i}$, and likewise the basis function $\phi_{i}^f$ associated with flux node $f_{i}$.

Eq. \ref{eq:conserv1d} is discretised in these function spaces pointwise at each of the solution nodes. To do this at each timestep:
\begin{enumerate}
    \item The piecewise degree $p$ solution polynomial is formed as $u(x) = \sum_i \phi_{i}^s(x) u(s_{i})$
    \item The flux $F$ is evaluated at each of the flux nodes $F_{i} = F\big( u(f_{i}) \big)$
    \item The piecewise degree $q$ flux polynomial is formed as $F(x) = \sum_i \phi_{i}^f(x) F_{i}$
    \item The derivative of the flux polynomial is evaluated at each solution node and used to update the solution
\end{enumerate}
In particular if we choose $p=q$ and the solution nodes to be equal to the flux nodes $\{ s_i \} = \{ f_i \}$, the procedure simplifies to an evaluation of the derivative of the flux function at each of the solution nodes.

While with this procedure the derivative of the flux polynomial $\frac{\partial F}{\partial x}$ is well defined in the interior of each element, it is in general multi-valued on the boundary of each element. Thus a unique value needs to be chosen for $\frac{\partial F}{\partial x}$ on element boundaries. For this model problem in Eq. \ref{eq:conserv1d} this is chosen simply as to be the value upwind to the velocity $a(u) = F'(u)$ at the boundary.

\begin{figure}
    \centering
    \includegraphics[scale=0.38]{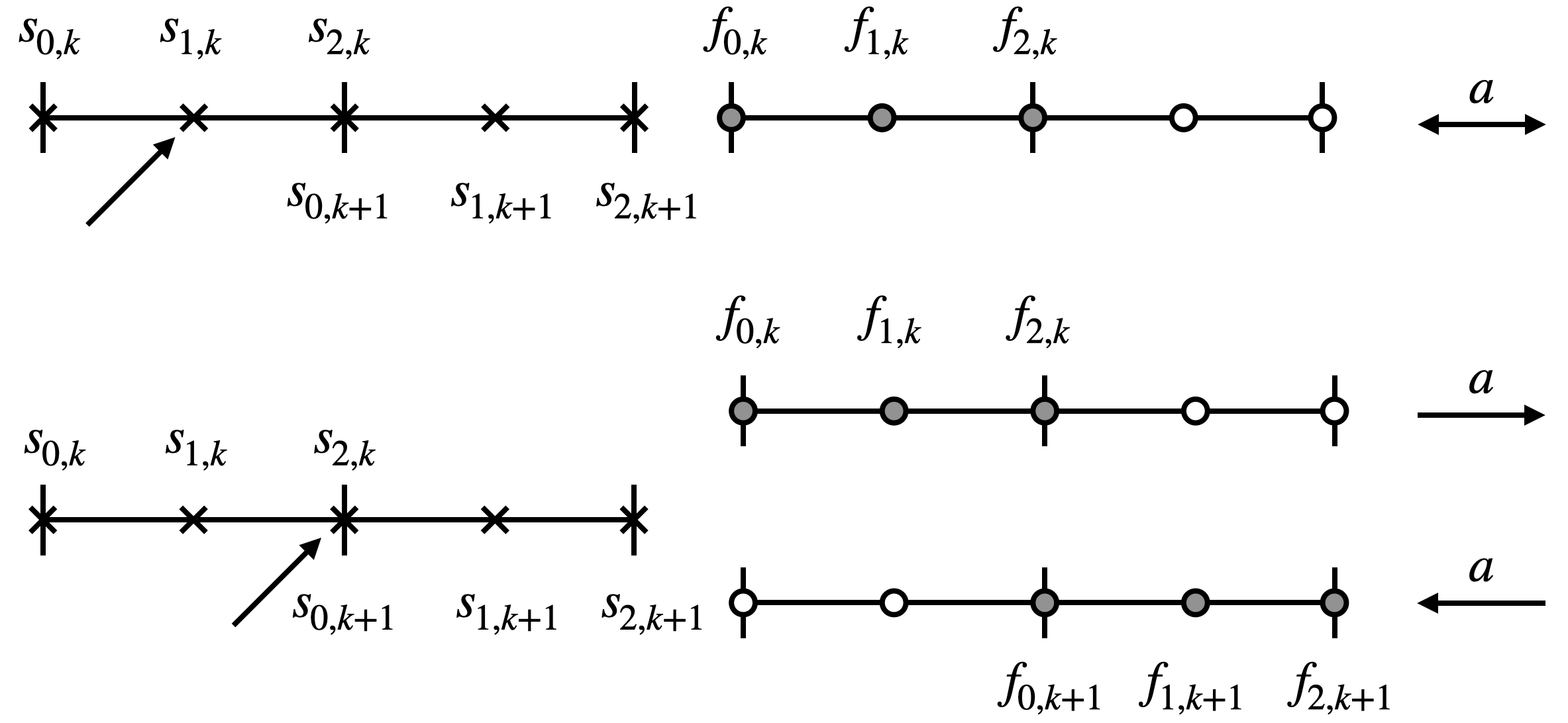}
    \caption{Schematic showing which shaded flux nodes (circles) are used to update the solution node (crosses) pointed to. On the top, solution nodes in the interior of an element are always updated using the flux nodes of that element regardless of the sign of the velocity $a$. On the bottom the flux nodes used to update the solution node on the boundary between the two elements is determined by the sign of the velocity $a$.}
    \label{fig:upwind1d}
\end{figure}

An example of this is shown in Fig. \ref{fig:upwind1d} for a $p=2$ mesh. For this example the solution nodes and flux nodes are both chosen to be equal $s_{i} = f_{i}$ and equidistant so that for all elements $K$
\begin{equation}
s_{i} - s_{i-1} = f_{i} - f_{i-1} = h
\end{equation}
Denoting the $i$-th solution node of the $k$-th element $K \in \mathcal{T}_h$ as $s_{i,k}$ and likewise the flux nodes $f_{i,k}$, in this case the derivative of the flux $F'$ at interior solution nodes $s_{1,k}$ is given simply as a spectral derivative using all the flux nodes in $K$
\begin{equation}
    F'(s_{1,k}) \approx \frac{1}{2h} F(f_{2,k}) - \frac{1}{2h} F(f_{0,k})
\end{equation}
For the derivative $F'$ at solution nodes on the boundary, the choice of stencil depends on the sign of the velocity $a(u) = F'(u)$. For instance at the solution node $s_{2,k}$ this is expressed as
\begin{equation} \label{eq:stencil1d_upwind}
    F'(s_{2,k}) \approx
    \begin{cases}
         \frac{3}{2h} F(f_{2,k}) - \frac{2}{h} F(f_{1,k}) + \frac{1}{2h} F(f_{0,k}), & \text{ } a(s_{2,k}) > 0 \\
         - \frac{3}{2h} F(f_{0,k+1}) + \frac{2}{h} F(f_{1,k+1}) - \frac{1}{2h} F(f_{2,k+1}),  & \text{ } a(s_{2,k}) \leq 0
    \end{cases}
\end{equation}
That is for a positive velocity to the right, the stencil from the upwind element from the left is used and vice versa. As can be seen from this example, in general the stencil for each point will contain at least $p+1$ points implying the method to be at least $p$-th order accurate in approximating the first derivative of an arbitrary smooth function.

\subsection{Connection to other methods}
\subsubsection{Spectral Differences}
It has been shown that many high-order methods such as DG, Nodal-DG, Spectral Volumes, Spectral Differences to name a few can all be unified under the Flux Reconstruction framework with different choices of reconstruction functions \cite{huynh2007flux}. It has been further shown in particular for 1D that the Spectral Volume and Spectral Difference methods are equivalent to one another as long as certain criteria are satisfied when constructing of each method \cite{wang2016spectral}.

We show that for constant-coefficient problems in 1D, that is when $a(u)$ is taken to be a constant, the FUSE method can be rewritten as a special case of the Spectral Difference (SD) method. This allows us to inherit properties including stability from Spectral Differences.

\begin{figure}
    \centering
    \includegraphics[scale=0.4]{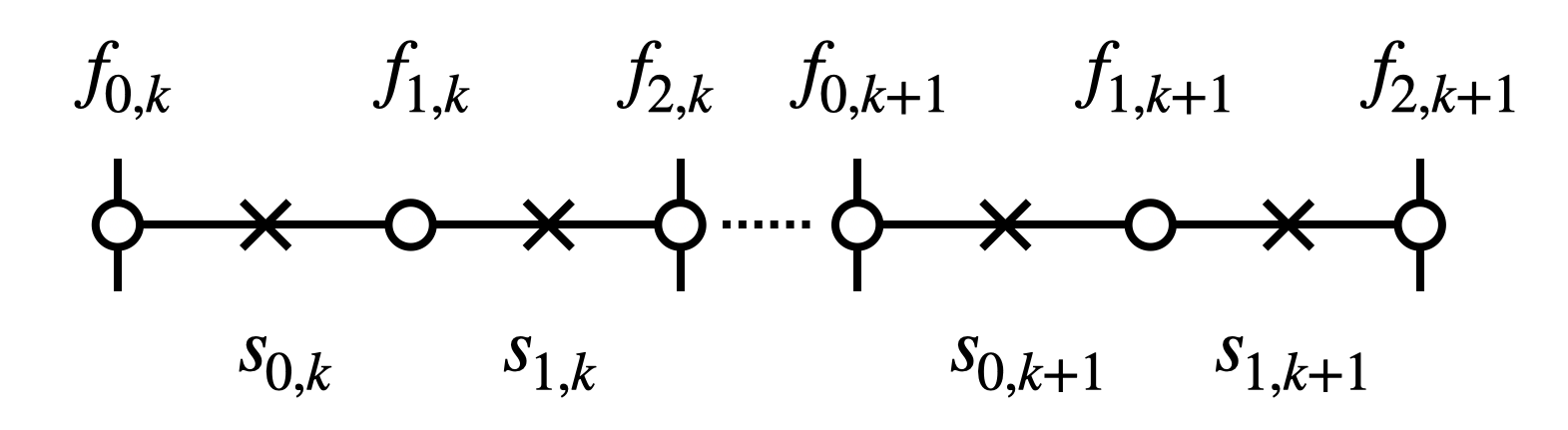}
    \caption{Spectral Difference method in 1D for $p=1$. Solution nodes are shown as before as crosses while flux nodes are shown as circles. Flux nodes are repeated on element boundaries.}
    \label{fig:sd}
\end{figure}

We briefly review the details of Spectral Differences in 1D for the model equation Eq. \ref{eq:conserv1d}. Similar to FUSE the domain is likewise divided up into distinct elements and a set of solution and flux nodes distributed within each element. Unlike FUSE in the SD method the number of flux nodes $q$ is fixed to always be one greater than the number of solution nodes $q = p+1$ so that the flux polynomial is of one degree higher than the solution polynomial. While there are no restrictions on the position of solution nodes, two flux nodes are always placed at each endpoint of the element. Unlike with our method however flux nodes at boundaries of elements are repeated resulting in a DG-like distribution of flux nodes. An example of the SD discretisation as described is shown in Fig. \ref{fig:sd}.

As with FUSE, to discretise Eq. \ref{eq:conserv1d} the SD method first reconstructs the solution $u$ and the flux $F(u)$ at each of the $p+1$ flux nodes. To determine the value of the flux at element boundaries a numerical flux function is applied. Finally the derivative flux function may then be calculated at each of the solution nodes via spectral differentiation. This construction was demonstrated to be linearly stable \cite{wang_sd08} for specific node distributions although stability was shown to depend only on the position of the flux nodes.

\begin{figure}
    \centering
    \includegraphics[scale=0.35]{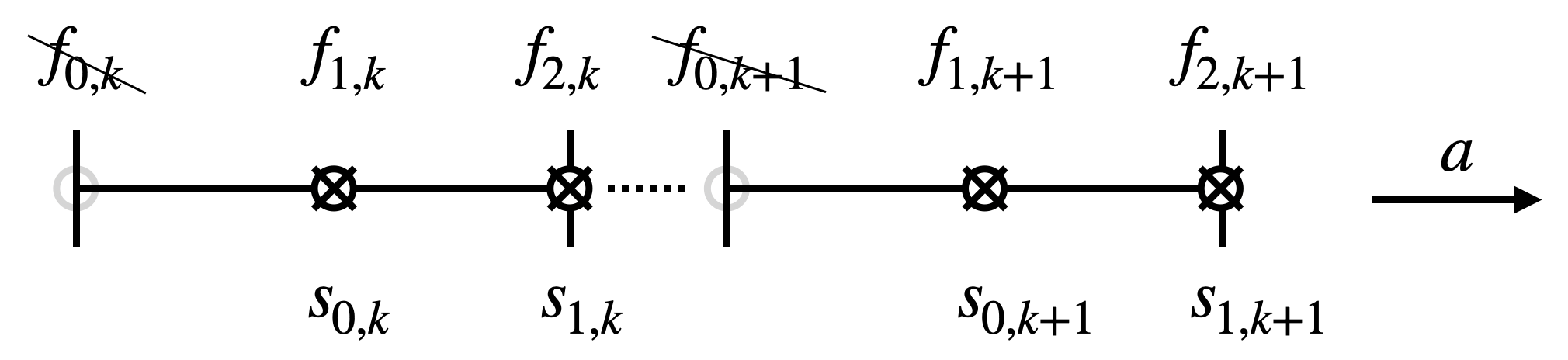}
    \caption{Schematic showing equivalence of our method to Spectral Differences in 1D for constant coefficient transport equation. Velocity $a(u) > 0$ is pointing to the right. Solution points are placed on top of flux points excluding one at left boundary in each element. Numerical flux is chosen to be the upwind flux so at each boundary the left value is taken and as such the leftmost flux value $f_{0,k}$ on each element is in effect ignored. The resulting Spectral Difference scheme is identical to our method in this specific case.}
    \label{fig:sd_equivalence}
\end{figure}

The equivalence of our method to Spectral Differences in this case can be seen in the following way and is shown in Fig. \ref{fig:sd_equivalence}. Assuming without loss of generality that the constant velocity $a(u) = 1$ the solution points $\{ s_0,...,s_{p} \}$ are placed on top of the flux points excluding the one of the leftmost boundary $\{ f_1,...,f_{p+1} \}$. This choice ensures that solution points do not overlap in space despite including an element boundary point. For the numerical flux at each boundary an upwind flux is be used, that is the value of the flux function from the element of the left is always taken. This choice of solution point location plus numerical flux means that the leftmost flux point on each element can be completely ignored, resulting in the flux points no longer in practice being repeated on element boundaries.

For a more general $a(u)$ the solution nodes cannot be placed in a way such that they do not overlap in space and for general flux functions repeated flux points on element boundaries cannot be ignored as was done with the above construction for this case. This means that in general the Spectral Difference solution approximated by discontinuous piecewise degree $p$ polynomials, the solution in our method is instead represented as a continuous piecewise degree $p$ polynomial analogous to Finite Element methods.

\subsubsection{Nodally integrated Petrov-Galerkin} \label{sect:petrov}
Many spectral element methods such as the Nodal-DG method can be understood as nodally-integrated Galerkin methods. In the case of Nodal-DG this is equivalent to a DG discretisation where the solution and flux nodes in each element are chosen to be equal to each other and to be Gauss-Lobatto nodes. Operators such as the mass and stiffness matrices are constructed via the usual Galerkin procedure, with the exception that numerical quadrature is computed using the Gauss-Lobatto solution nodes themselves. This results in an inexact diagonal mass matrix that is much simpler to invert than that of full DG.

In the same spirit, the FUSE method can be reformulated as an nodally-integrated Petrov-Galerkin method. Given a mesh $\mathcal{T}_h$ equipped with a Finite Element space $V(\mathcal{T}_h)$, to discretise Eq. \ref{eq:conserv1d} we multiply both sides with a set of test functions $\{\psi_i\} \in V$ and integrate over the domain
\begin{equation} \label{eq:1d_weak}
    \int_\Omega \frac{\partial u}{\partial t} \psi_i + \frac{\partial F(u)}{\partial x} \psi_i dx = 0
\end{equation}
Substituting in the form of the solution $u = \sum_j u(s_{j}) \phi_{j}^s$ and the flux $F = \sum_l F(f_{l}) \phi_{l}^f$ gives the linear system
\begin{equation} \label{eq:1d_weak_basis}
     \sum_j \int_\Omega \psi_i \phi_j^s dx \cdot \frac{\partial}{\partial t}u_j \ + \sum_l \int_\Omega \psi_i \frac{d \phi_l^f}{dx} dx \cdot F(f_l) = 0
\end{equation}

\begin{figure}
    \centering
    \includegraphics[scale=0.4]{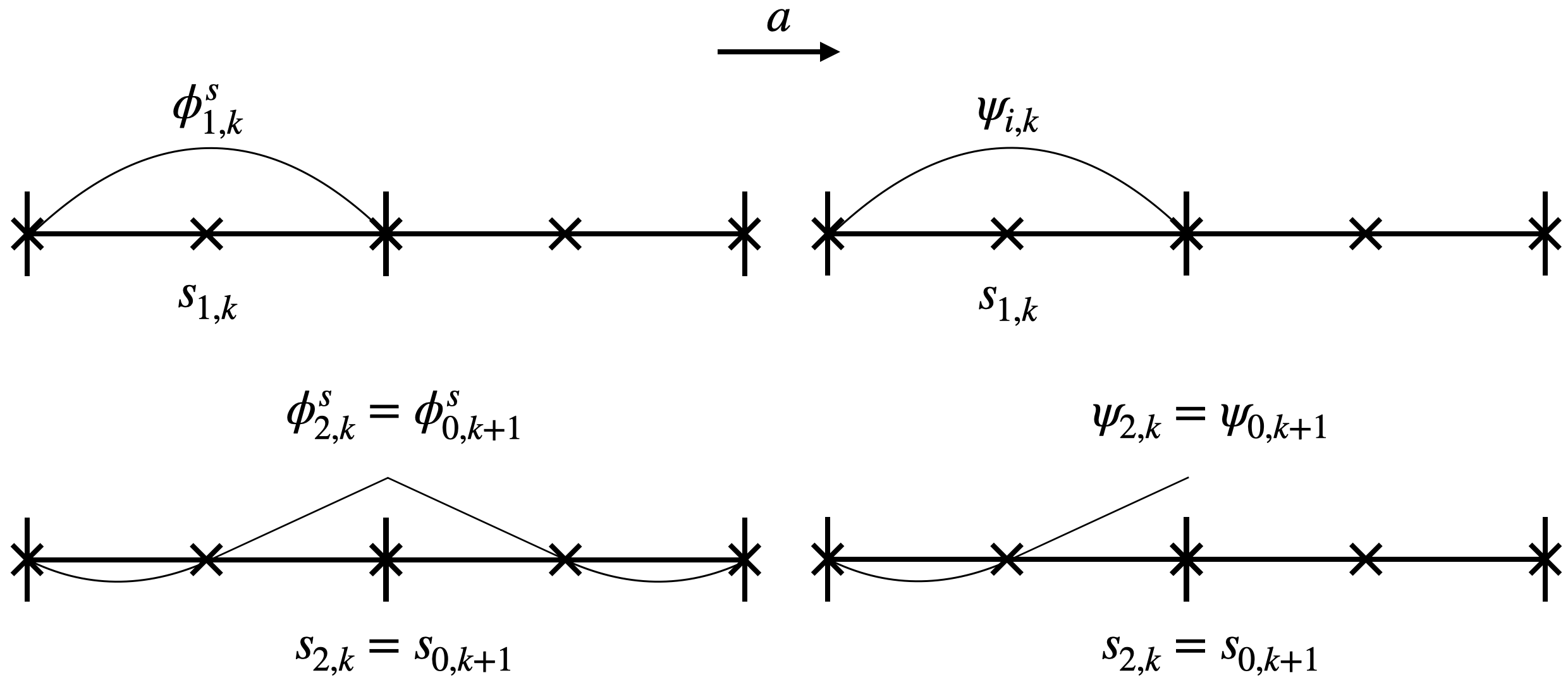}
    \caption{Schematic of FUSE Petrov-Galerkin test functions on two element $p=2$ mesh. Velocity $a$ is assumed to be positive. Top row shows basis and test functions for interior solution node $s_{1,k}$: on left its basis function $\phi_{1,k}^s$ on right its test function $\psi_{1,k} = \phi_{1,k}^s$. Bottom row shows basis and test functions for solution node $s_{i+1}$: on left its basis function $\phi_{2,k}^s$ which has support on both left and right elements, on right its test function $\psi_{2,k}$ which is its basis function with support restricted to the left upwind element.}
    \label{fig:upwind_fn}
\end{figure}

In a standard Galerkin method, the test functions are chosen to be the same basis functions used to discretise the solution. For the FUSE construction however these test functions $\psi_i$ are instead ``upwinded'' basis functions. For solution nodes inside an element, its basis function has support only on that element and so the test functions for these nodes are exactly the basis functions themselves. However basis functions for solution nodes on the boundary elements have support on all elements that they border. In this case the upwinded test function is simply the solution basis function restricted to the element upwind from the velocity. A schematic of this is shown in Fig. \ref{fig:upwind_fn}. 

Having chosen test functions $\{ \psi_i \}$, Eq. \ref{eq:1d_weak_basis} is then constructed via nodal integration at the solution points themselves. This is  in contrast to standard FEM where usually a new separate set of quadrature points are introduced. Under this integration rule Eq. \ref{eq:1d_weak_basis} becomes
\begin{equation}
     \sum_l \sum_j \psi_i(s_j) \phi_j^s(s_l) w_l \cdot \frac{\partial}{\partial t}u_j \ + \sum_l \sum_k \psi_i (s_l) \frac{d \phi_k^f}{dx} (s_l) w_l \cdot F(u_k) = 0
\end{equation}
where $w_l$ is the integration weight associated with solution node $s_i$. This can be simplified using the following two facts: 1) $\phi_j^s(s_l) = \delta_{jl}$, 2) $\psi_i = \phi_i^s$ on the support of $\psi_i$, to give
\begin{equation}
      w_i \bigg( \frac{\partial}{\partial t}u_i \ + \sum_k \mathbbm{1}_{\text{supp}(\psi_i)} \frac{d \phi_k^f}{dx} (s_i) \cdot F(u_k) \bigg) = 0
\end{equation}
Dividing this equation by the integration weight $w_i$ we recover the FUSE method. 

The FUSE method shares similarities also with stabilised continuous Galerkin methods such as SUPG and VMS, with the key one being that all these methods use upwinded test functions for stabilisation. However a difference lies in that no explicit user-defined upwind parameter is required in the formulation of the FUSE method that is present for instance in SUPG. Furthermore as we show in the next section, while stability of methods such as SUPG are independent of solution node locations, much like Nodal-DG stability of FUSE is highly dependent on the position of the solution and flux nodes.

\subsection{Stability} \label{sec:stability1d_convect}
We perform an analysis similar to that performed for DG in \cite{hu1999analysis} and Spectral Differences in \cite{wang_sd08} to show stability of our method. For this we assume without loss of generality that the velocity in Eq. \ref{eq:non_conserv_convect} that  is a positive constant $a(u) = 1$ and consider the linear advection equation
\begin{equation} \label{eq:transport1d}
    \frac{\partial u}{\partial t} + \frac{\partial u}{\partial x} = 0
\end{equation}
Denoting the discrete solution on the domain $\Omega$ as $u_h$, the above equation can then be written as
\begin{equation}
    \frac{\partial u_h}{\partial t} + A u_h = 0
\end{equation}
where $A$ is the discrete upwind first derivative operator formed from the procedure in the previous subsections. Applying an inner product with $u_h$ to all terms in the equation we get
\begin{equation}
    \frac{\partial \| u_h \|_2^2}{\partial t} + 2u_h^T A u_h = 0
\end{equation}
using the definition of the operator norm of $A$ this is equivalent to
\begin{equation}
    \frac{\partial \| u_h \|_2^2}{\partial t} + 2 \| A \|_2 \| u_h \|_2^2 =0
\end{equation}
and applying Gronwall's inequality we find that
\begin{equation}
    \| u_h \|_2^2 \leq C e^{ -2\| A \|_2 t }
\end{equation}
Thus for stability of this discretisation we require all eigenvalues of the operator $A$ to have strictly non-negative real part.

We observe in numerical experiments the eigenvalues of the linear operator $A$ to be independent on the choice of flux nodes and depend only on solutions nodes. This is in contrast to SD, where stability of the linear operator depends only on the position of flux nodes. This behaviour can be explained as follows: under the action of the operator $A$ each solution node is updated using only the $p+1$ solution nodes and its associated basis functions of one single element. As we pick in the FUSE framework the number of flux points to be greater or equal to the number of solution points $q \geq p$, this implies that $\mathcal{V}_s \subseteq \mathcal{V}_f$, and so the flux polynomial space on each element contains the solution space. However as each element consists of $p+1$ nodes with associated polynomial basis functions $\{ \phi_i^s \}$ of degree $p$, these $\{ \phi_i^s \}$ already uniquely define a complete polynomial basis for $\mathcal{V}_s$. Thus modifying the flux nodes has no effect on $\mathcal{V}_s$ on an element and thus does not affect the stability of the linear operator $A$.

As a result for simplicity we simply choose the number of solution nodes to be equal to the number of flux nodes $p=q$, and for the sets of nodes to be the same $\{ s_i \} = \{ f_i \}$. For the remainder of the paper no distinction is as a result made between the two sets of nodes.

\begin{figure}
    \centering
    \includegraphics[scale=0.4]{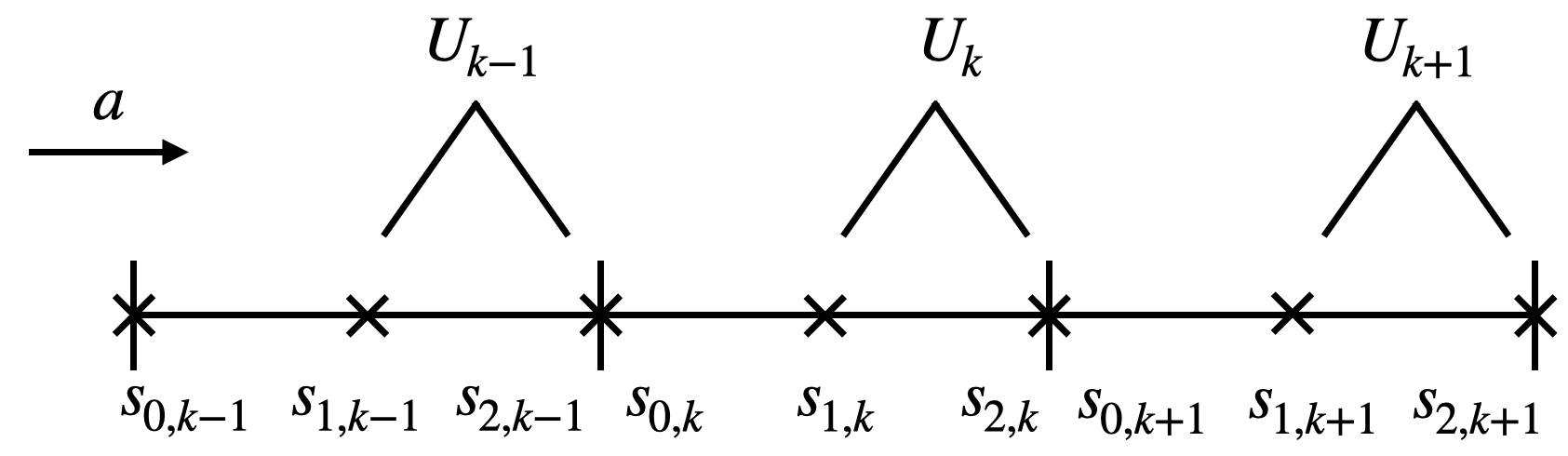}
    \caption{Setup for stability analysis of upwind discretisation for three $p=2$ elements. Solution nodes shown here are indexed as $s_{j,k}$ for the $j$-th node of the $k$-th element. Note that solution node $2$ of element $k$ is the same as solution node $0$ of element $k+1$. With the velocity assumed to be positive, solution nodes $s_{1,k}$ and $s_{2,k}$ of each element are grouped together into a vector $U_{k}$ for the Von Neumann analysis.}
    \label{fig:vonneumann_1d}
\end{figure}

We inspect the eigenvalues of the upwind operators $A$ by adapting an approach employed usually in von Neumann analysis for Finite Difference methods, whereby the solution is taken to be a linear combination of eigenfunctions of the first derivative operator $W(x) = e^{i x\xi}$ with $\xi \in [0,2\pi]$ the wavenumber. We first label the solution nodes of the $k$-th element as $\{ s_{i,k} ~, i = 0,1,...,p \},$ and group up the solution values $u_h$ on the nodes $\{ s_{1,k}, ... , s_{p,k} \}$ as a vector $U_k = \bigg( u(s_{1,k}),...,u(s_{p,k}) \bigg)^T$ as shown in Fig. \ref{fig:vonneumann_1d}. As the velocity here is assumed positive, every solution node in the $k$-th element has a stencil that only involves other solution nodes in that element. This implies that

\begin{equation}
    A U_k = \bar{A}U_k + A_1 U_{k-1}
\end{equation}
where $A_1$ is a rank one matrix with non-zero entries only in the last column. To illustrate this we refer again to the $p=2$ example in Fig. \ref{fig:upwind1d} and Eq. \ref{eq:stencil1d_upwind}. For this example the linear system becomes
\begin{equation}
    A U_k
    =
    \begin{pmatrix}
        \frac{3}{2h} & -\frac{2}{h} \\
        \frac{1}{2h} & 0
    \end{pmatrix}
    U_k
    +
    \begin{pmatrix}
        0 & \frac{1}{2h} \\
        0 & -\frac{1}{2h}
    \end{pmatrix}
    U_{k-1}
\end{equation}
Using the assumption that the solution is of the form $W(x) = e^{ix\xi}$ as stated above, this gives that
\begin{equation}
    A U_k
    =
    \Bigg[
    \begin{pmatrix}
        \frac{3}{2h} & -\frac{2}{h} \\
        \frac{1}{2h} & 0
    \end{pmatrix}
    +
    e^{-i2h\xi}
    \begin{pmatrix}
        0 & \frac{1}{2h} \\
        0 & -\frac{1}{2h}
    \end{pmatrix}
    \Bigg]
    U_{k}
\end{equation}
which allows us to consider the spectrum of $A$ as a function of $\xi$ the wavenumber.

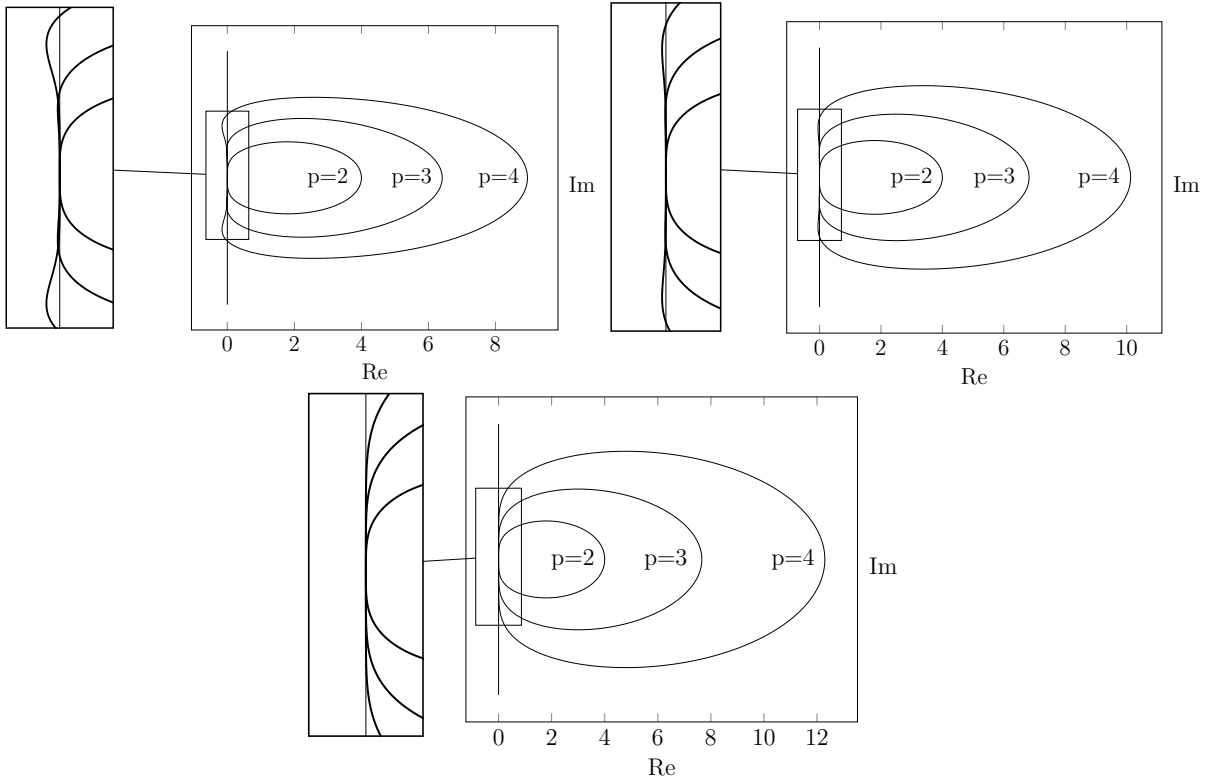
\begin{figure}[h]
    \centering
    \begin{minipage}{0.48\textwidth}
    \resizebox{\linewidth}{!}{%
    \begin{tikzpicture}[font=\large]
    \begin{scope}[
        spy using outlines={
            rectangle,
            magnification=2.5,
            connect spies,
            width=2cm,height=6cm,
        },
        ]
        \begin{axis}
        [
            xlabel = Re,
            ylabel = Im,
            ytick = \empty,
            y label style={at={(axis description cs:1.12,0.48)},rotate=270},
        ]
        
        \addplot[mark=none, color=black] 
        table[ x=real_uniform_p2, y=imag_uniform_p2 ]{vonneumann_uniform_p2.txt};

        \addplot[mark=none, color=black] 
        table[ x=real_uniform_p3, y=imag_uniform_p3 ]{vonneumann_uniform_p3.txt};

        \addplot[mark=none, color=black] 
        table[ x=real_uniform_p4, y=imag_uniform_p4 ]{vonneumann_uniform_p4.txt};

        \addplot[mark=none,domain=0:4,very thin] coordinates {(0, -10) (0, 10)};

        \node[] at (axis cs: 3.0,0.0) {p=2};
        \node[] at (axis cs: 5.5,0.0) {p=3};
        \node[] at (axis cs: 8.1,0.0) {p=4};

        \coordinate (spy point) at (axis cs:-5.0,0.8);
        \coordinate (point) at (axis cs:0.0,0.2);
        \spy on (point) in node (spy) at (spy point);
        \end{axis}
        \end{scope}
    \end{tikzpicture}%
    }
    \end{minipage}
    \begin{minipage}{0.48\textwidth}
    \resizebox{\linewidth}{!}{%
    \begin{tikzpicture}[font=\large]
    \begin{scope}[
        spy using outlines={
            rectangle,
            magnification=2.5,
            connect spies,
            width=2cm,height=6cm,
        },
        ]
        \begin{axis}
        [
            xlabel = Re,
            ylabel = Im,
            ytick = \empty,
            y label style={at={(axis description cs:1.12,0.48)},rotate=270},
        ]
        
        \addplot[mark=none, color=black] 
        table[ x=real_uniform_p2, y=imag_uniform_p2 ]{vonneumann_uniform_p2.txt};

        \addplot[mark=none, color=black] 
        table[ x=real_gausslobatto_p3, y=imag_gausslobatto_p3 ]{vonneumann_gll_p3.txt};

        \addplot[mark=none, color=black] 
        table[ x=real_gausslobatto_p4, y=imag_gausslobatto_p4 ]{vonneumann_gll_p4.txt};

        \addplot[mark=none,domain=0:4,very thin] coordinates {(0, -10) (0, 10)};

        \node[] at (axis cs: 3.0,0.0) {p=2};
        \node[] at (axis cs: 5.7,0.0) {p=3};
        \node[] at (axis cs: 9.1,0.0) {p=4};

        \coordinate (spy point) at (axis cs:-5.0,0.8);
        \coordinate (point) at (axis cs:0.0,0.2);
        \spy on (point) in node (spy) at (spy point);
        \end{axis}
        \end{scope}
    \end{tikzpicture}%
    }
    \end{minipage}
    \begin{minipage}{0.48\textwidth}
    \resizebox{\linewidth}{!}{%
    \centering
    \begin{tikzpicture}[font=\large]
    \begin{scope}[
        spy using outlines={
            rectangle,
            magnification=2.5,
            connect spies,
            width=2cm,height=6cm,
        },
        ]
        \begin{axis}
        [
            xlabel = Re,
            ylabel = Im,
            ytick = \empty,
            y label style={at={(axis description cs:1.12,0.48)},rotate=270},
        ]
        
        \addplot[mark=none, color=black] 
        table[ x=real_uniform_p2, y=imag_uniform_p2 ]{vonneumann_uniform_p2.txt};

        \addplot[mark=none, color=black] 
        table[ x=real_legendre_p3, y=imag_legendre_p3 ]{vonneumann_gl_p3.txt};

        \addplot[mark=none, color=black] 
        table[ x=real_legendre_p4, y=imag_legendre_p4 ]{vonneumann_gl_p4.txt};

        \addplot[mark=none,domain=0:4,very thin] coordinates {(0, -10) (0, 10)};

        \node[] at (axis cs: 2.8,0.0) {p=2};
        \node[] at (axis cs: 6.3,0.0) {p=3};
        \node[] at (axis cs: 11.1,0.0) {p=4};

        \coordinate (spy point) at (axis cs:-5.0,-0.4);
        \coordinate (point) at (axis cs:0.0,0.2);
        \spy on (point) in node (spy) at (spy point);
        \end{axis}
        \end{scope}
    \end{tikzpicture}%
    }
    \end{minipage}
    \caption{Comparison of spectrum of 1st order upwind operator for different node distributions. Here the solution and flux nodes are chosen to be same. The uniform (top left) and Gauss-Lobatto (top right) nodes display a small instability whilst the Gauss-Legendre nodes plus endpoints (bottom) are shown to be stable.}
    \label{fig:eigenvalues_1d}
\end{figure}
The eigenspectra for the upwind first derivative differential operator using three different solution node distributions for $h=1$ are shown in Fig. \ref{fig:eigenvalues_1d}. For $p>2$, we observe that the eigenspectra is not contained in the positive half plane using a uniform node distribution resulting in an unstable operator. Furthermore in contrast to Nodal-DG, we observe the first derivative operator using Gauss-Lobatto nodes to also be unstable.

We instead consider the set of nodes introduced by Van den Abeele \cite{wang_sd08} for Spectral Differences, defined on the reference domain $[-1,1]$ as the standard Gauss-Legendre nodes used in Gaussian quadrature plus the two endpoints at $\{-1,1\}$. These nodes for $p=2,3,4$ are shown in Table \ref{table:nodes}. The stability of the Spectral Difference method using this distribution of flux nodes has been proven by Jameson \cite{jameson2009}. For the FUSE method, we have verified numerically the stability of the upwind first derivative operator with this distribution of solution nodes for orders $p \leq 20$.
\begin{table}
\begin{center}
    \begin{tabular}{ |c|c|c|c| } 
    \hline
    $p$& 2 & 3 & 4 \\
    \hline
    & & & \\
    $\{ s_i \}$ & -1, 0, 1 & -1, $-\sqrt{\frac{1}{3}}, \sqrt{\frac{1}{3}}$, 1 & -1, $-\sqrt{\frac{3}{5}},0,\sqrt{\frac{3}{5}},$ 1 \\ [15pt]
    \hline
    \end{tabular}
    \caption{Table of Gauss-Legendre nodes plus endpoints for $p=2,3,4$.}
    \label{table:nodes}
\end{center}
\end{table}

\subsection{Conservation} 
We present a method to demonstrate conservation of the FUSE method for $p \geq 3$ that borrows from the framework of Finite Volume methods (FVM) which are well-known to be conservative. Integrating Eq. \ref{eq:conserv1d} in space over any arbitrary interval $[a,b]$
\begin{equation}
    \frac{1}{|b-a|} \frac{\partial}{\partial t} \int_a^b u ~dx + \frac{1}{|b-a|}\int_a^b \frac{\partial F(u)}{\partial x} ~dx = 0
\end{equation}
using the fundamental theorem of calculus this can be simplified to
\begin{equation} \label{eq:fv1d}
    \frac{\partial \bar{u}}{\partial t} + \frac{1}{|b-a|} \bigg( F\big(u(b)\big) - F\big(u(a)\big) \bigg) = 0
\end{equation}
with the volume averaged solution over the interval introduced as $\bar{u} = \frac{1}{|b-a|} \int_a^b u ~dx$. For standard FVM the domain is partitioned into distinct discrete volume cells each with a defined cell average, and the solution propagated forward in time using Eq. \ref{eq:fv1d} along with a time integrator of choice. Conservation of the cell averages $\bar{u}$ is acheived discretely as Eq. \ref{eq:fv1d} ensures that any flux exiting a cell boundary is identical to the one entering the adjacent cell.

\begin{figure}
    \centering
    \includegraphics[scale=0.38]{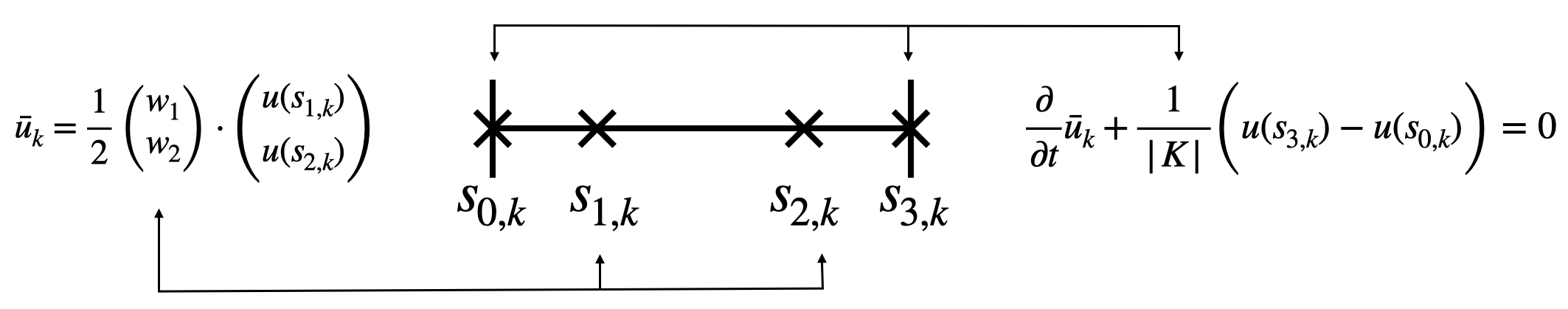}
    \caption{Schematic of FVM cell average construction for a single $p=3$ element $K$. The cell averages are calculated using values at the interior solution nodes which are the standard Gauss-Legendre points using Gaussian quadrature. The endpoints are used to update the cell averages in time.}
    \label{fig:fv1d}
\end{figure}

We consider the following FVM inspired construction on top of our method: each element in $\mathcal{T}_h$ is considered to be a volume cell with a defined cell average. For the $k$-th element $K \in \mathcal{T}_h$ this is defined as $\bar{u}_k = \frac{1}{|K|}\int_{K} u ~dx$. As the solution in each cell can be written as a degree $p$ polynomial, these cell averages can be calculated exactly using the interior solution nodes in each element as they are chosen to be the Gauss-Legendre quadrature points 
\begin{equation} \label{eq:conserv_constraint}
    \bar{u}_k = \sum_{i=1}^{p-1} u(s_{i,k}) \cdot \frac{w_i}{2}
\end{equation}
where $w_i$ are the standard Gauss-Legendre quadrature weights and $s_{i,k}$ as before the solution nodes of element $K$. At each timestep the solution nodes are updated as per usual with the first order operators previously described. The cell averages are on the other hand updated using the solution values at the endpoints of each element with Eq. \ref{eq:fv1d}. This is shown in Fig. \ref{fig:fv1d}.

The goal is to show that the equations for updating cell average can in fact be written as a linear combination of the solution node updates. To show this we can focus solely on the $k$-th element element $K$. On the interior of the element for solution nodes the differential equation Eq. \ref{eq:transport1d} is satisfied pointwise exactly for a degree $p$ polynomial. A linear combination of these equations can be taken as follows using the corresponding Gauss integration weight
\begin{equation}
    \sum_{i=1}^p \frac{w_i}{2} \bigg( \frac{\partial u(s_{i,k})}{\partial t} + \frac{\partial F\big(u(s_{i,k})\big)}{\partial x} \bigg) = 0
\end{equation}
which on the element $K$ is the discrete analogue of the integral equation
\begin{equation}
    \int_{K} \frac{\partial u}{\partial t} + \frac{\partial F(u)}{\partial x} ~dx = 0
\end{equation}
which can be rewritten as
\begin{equation}
    \frac{\partial \bar{u}_k}{\partial t} + \frac{1}{|K|} \bigg( F\big(u(s_{p,k})\big) - F\big(u(s_{0,k})\big) \bigg) = 0
\end{equation}
and summing over the entire domain and utilising the periodic boundary conditions
\begin{equation}
    \sum_{k} \frac{\partial \bar{u}_k}{\partial t} = 0
\end{equation}
implying that average over all cells is conserved. As the solution nodes are updated independently of the cell averages, this implies that the original method implicitly conserves cell averages. In practice the FVM-type construction need not be explicitly constructed as the cell averages are automatically conserved and can easily be obtained at any timestep for each element $K$ via Eq. \ref{eq:conserv_constraint}.

\section{2nd derivative operators in 1D} \label{sect:poi1d}
\subsection{Preliminaries}
We consider the following model problem for the construction of second derivative operators
\begin{equation} \label{eq:heat1d}
    \frac{\partial u}{\partial t} - \Delta u = 0
\end{equation}
on the domain $\Omega = [0,1]$ with periodic boundary conditions. Given that our discretisation results in nodes consistent with FEM possible choices for discretising the Laplace operator could be with continuous Finite Element or Spectral Element methods. In this paper we also present an alternative method for doing so inspired by Local Discontinuous Galerkin (LDG) operator splitting \cite{ldg1998} that is consistent with upwind first derivative operators.

\subsection{Operator splitting}
\begin{figure}
    \centering
    \includegraphics[scale=0.35]{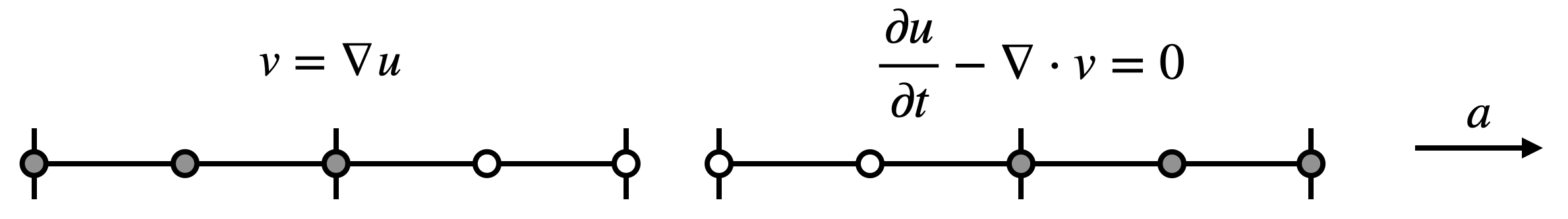}
    \caption{An example LDG type operator splitting for second derivative equations. The velocity $a$ is arbitrarily picked to be positive. To update the indicated solution node, Eq. \ref{eq:ldg1} is discretised using an upwind operator, so nodes on the left element are used. Eq. \ref{eq:ldg2} the didscretised using the downwind operator, so nodes on right element are used.}
    \label{fig:ldg1d}
\end{figure}

To discretise the Laplace operator a new variable is introduced for the gradient $v=\nabla u$ and Eq. \ref{eq:heat1d} rewritten as
\begin{align}
    \label{eq:ldg1} v &= \nabla u \\ 
    \label{eq:ldg2} \frac{\partial u}{\partial t} - &\nabla \cdot v = 0 
\end{align}
This is a system of differential equations which can be discretised using the procedure described in the previous section. To do so a velocity is arbitrarily prescribed at each point of the mesh, and the operator in first equation Eq. \ref{eq:ldg1} constructed to be upwind to the prescribed velcocity, and the second equation Eq. \ref{eq:ldg1} downwind to the velocity. An example of this is shown in Fig. \ref{fig:ldg1d}. Discretising as in this upwind-downwind fashion results in a linear system of the form
\begin{align}
    v_h &= A^{+}u_h \\
    \frac{\partial u_h}{\partial t} - &A^{-}v_h = 0
\end{align}
where $u_h, v_h$ denotes the discrete solution and gradient, and $A^+, A^-$ the upwind and downwind first derivative operators respectively. This then gives the overall second derivative operator $A$ as
\begin{equation}
    \frac{\partial u_h}{\partial t} - Au_h = 0, \hspace{10mm}  A = A^{-}A^{+}
\end{equation}

\subsection{Stability}
The stability of the resulting discrete Laplace operator is analysed in the same manner as for first derivative operators above it is assumed that the solution $u$ can be written as a linear combination of functions of the form $W(x) = e^{ix\xi}$. Once again we we consider the example $p=2$ in Fig. \ref{fig:upwind1d} where the solution nodes are assumed to be equal to the flux nodes and the velocity assumed positive everywhere. In this case the first derivatives of the solution $u$ and the gradient $v$ at solution nodes $\{ s_{j,k} \}$ of an element $K$ can be calculated as
\begin{align}
    \frac{\partial u(s_{1,k})}{\partial x} &= \frac{1}{2h} u(s_{2,k}) - \frac{1}{2h}u(s_{2,k-1}) \\
    \frac{\partial u(s_{2,k})}{\partial x} &= \frac{3}{2h} u(s_{2,k}) -\frac{2}{h}u(s_{1,k}) + \frac{1}{2h}u(s_{2,k-1}) \\
    \frac{\partial v(s_{1,k})}{\partial x} &= \frac{1}{2h} v(s_{2,k}) - \frac{1}{2h}v(s_{2,k-1}) \\
    \frac{\partial v(s_{2,k})}{\partial x} &= -\frac{3}{2h} v(s_{2,k}) +\frac{2}{h}v(s_{1,k+1}) - \frac{1}{2h}v(s_{2,k+1})
\end{align}
which can be written as the following linear system
\begin{align}
    A^{+}U_k &=
    \begin{pmatrix}
        0 & \frac{1}{2h} \\
        -\frac{2}{h} & \frac{3}{2h}
    \end{pmatrix}
    U_k
    +
    \begin{pmatrix}
        0 & -\frac{1}{2h} \\
        0 & \frac{1}{2h}
    \end{pmatrix}
    U_{k-1} \\
    A^{-}V_k &=
    \begin{pmatrix}
        0 & \frac{1}{2h} \\
        0 & -\frac{3}{2h}
    \end{pmatrix}
    V_k
    +
    \begin{pmatrix}
        0 & -\frac{1}{2h} \\
        0 & 0
    \end{pmatrix}
    V_{k-1}
    +
    \begin{pmatrix}
        0 & 0 \\
        \frac{2}{h} & -\frac{1}{2h}
    \end{pmatrix}
    V_{k+1}
\end{align}
combining the two gives the overall system
\begin{equation}
    AU_k =
    \begin{pmatrix}
    0& -\frac{1}{4h^2} \\
    0& 0
    \end{pmatrix}
    U_{k-2}
    +
    \begin{pmatrix}
    \frac{1}{h^2}& -\frac{1}{2h^2} \\
    0& -\frac{3}{4h^2}
    \end{pmatrix}
    U_{k-1}
    +
    \begin{pmatrix}
    -\frac{1}{h^2}& \frac{3}{4h^2} \\
    \frac{3}{h^2}& -\frac{7}{2h^2}
    \end{pmatrix}
    U_{k}
    +
    \begin{pmatrix}
    0& 0 \\
    \frac{1}{h^2}& \frac{1}{4h^2}
    \end{pmatrix}
    U_{k+1}
\end{equation}
again using the assumption that the solution is of the form $W(x) = e^{ix\xi}$ as stated above, this gives that
\begin{equation}
    A U_k
    =
    \Bigg[e^{-i4h\xi}
    \begin{pmatrix}
    0& -\frac{1}{4h^2} \\
    0& 0
    \end{pmatrix}
    +
    e^{-i2h\xi}
    \begin{pmatrix}
    \frac{1}{h^2}& -\frac{1}{2h^2} \\
    0& -\frac{3}{4h^2}
    \end{pmatrix}
    + \begin{pmatrix}
    -\frac{1}{h^2}& \frac{3}{4h^2} \\
    \frac{3}{h^2}& -\frac{7}{2h^2}
    \end{pmatrix}
    + 
    e^{i2h\xi}
    \begin{pmatrix}
    0& 0 \\
    \frac{1}{h^2}& \frac{1}{4h^2}
    \end{pmatrix}
    \Bigg]
    U_{k}
\end{equation}

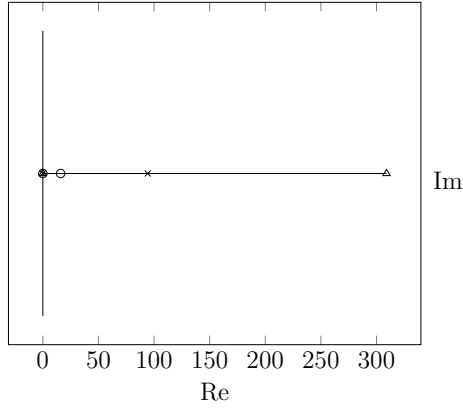
\begin{figure}
    \centering
    \begin{tikzpicture}[scale=0.8,font=\large]
        \begin{axis}
        [
            xlabel = Re,
            ylabel = Im,
            ytick = \empty,
            y label style={at={(axis description cs:1.12,0.48)},rotate=270},
        ]
        
        \addplot[mark=o, color=black] 
        table[ x=x, y=y ]{poi_p2.txt};
        \label{plot:vnpoi2}

        \addplot[mark=x, color=black] 
        table[ x=x, y=y ]{poi_p3.txt};
        \label{plot:vnpoi3}

        \addplot[mark=triangle, color=black] 
        table[ x=x, y=y ]{poi_p4.txt};
        \label{plot:vnpoi4}

        \addplot[mark=none,domain=0:4,very thin] coordinates {(0, -1) (0, 1)};

        \end{axis}
    \end{tikzpicture}%
    \caption{Spectrum of Laplace order upwind operator with Gauss-Legendre points plus endpoints node distribution for solution and flux nodes. \ref{plot:vnpoi2} shows distribution for $p=2$, \ref{plot:vnpoi3} shows distribution for $p=3$, \ref{plot:vnpoi4} shows distribution for $p=4$.}
    \label{fig:vonneumann_poisson}
\end{figure}
The eigenspectra for several degrees $p$ is shown in Fig. \ref{fig:vonneumann_poisson} for $h=1$. For these plots the solution and flux points are equal have been chosen as the Gauss-Legendre nodes plus endpoints. Stability of the upwind-downwind Laplace operator is established as all the eigenvalues lie in the positive half-plane.

\section{Derivative operators in higher dimensions}
\subsection{Preliminaries}
To discretise any function $u$ in higher dimensions, the domain is split into distinct elements on which a set of nodes is distributed. For the purposes of this paper we focus only on $2^d$ sided elements with $d$ the dimension (i.e. quadrilaterals in 2D, hexahedra in 3D) allowing us to define nodes in higher dimension simply as an outer product of a given one dimensional node set. Specifically in $d$-dimensions and given a polynomial order $p$, we pick the nodeset to be an $d$-dimensional outer product of the $p+1$ Gauss-Legendre nodes plus endpoints defined in Sect. \ref{sec:stability1d_convect}. As with the 1D case nodes on element boundaries are not repeated and are shared between neighbouring elements resulting in a nodes consistent with FEM. Furthermore as mentioned in Sect. \ref{sec:stability1d_convect} no distinctions are made between solution and flux nodes as they are chosen to be equal to one another.

\subsection{1st derivative operators} \label{sect:2d_1stderiv}
For higher dimensions the general first order equation in conservative form can be written as
\begin{equation} \label{eq:conv_2d}
    \frac{\partial u}{\partial t} + \nabla \cdot F(u) = 0
\end{equation}
on some domain $\Omega \subset \mathbb{R}^d$ where $d$ denotes the dimension. Assuming sufficient continuity on $f$ this can be rewritten using the chain rule
\begin{equation}
    \frac{\partial u}{\partial t} + a(u) \cdot \nabla u = 0
\end{equation}
where $a(u) = F'(u)$. As with the 1D case $a(u)$ can be interpreted as the velocity used pointwise to determine the upwind direction.

\begin{figure}
    \centering
    \includegraphics[scale=0.35]{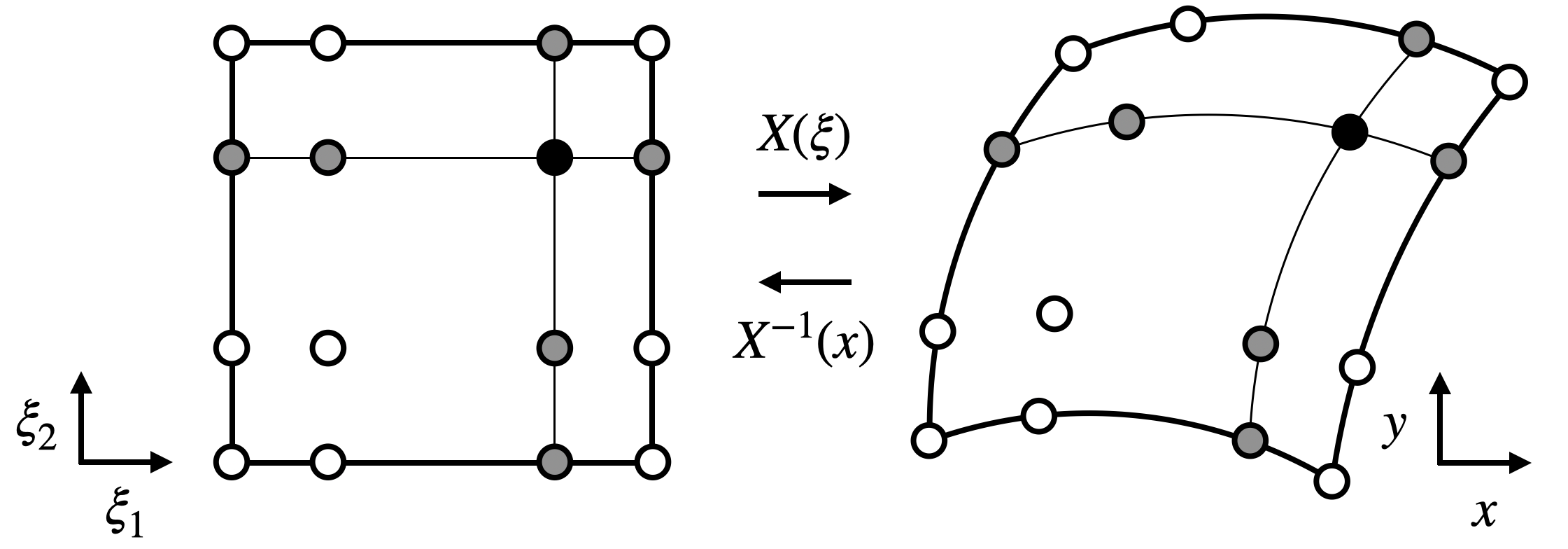}
    \caption{Stencil construction for points within an element for $p=3$. Elements are mapped from physical space to reference space $\boldsymbol{\xi} = \boldsymbol{X}^{-1}(\boldsymbol{x})$ where stencils are constructed using only the two 1D lines defined by the tensor product node structure. }
    \label{fig:stencils2d}
\end{figure}

To form derivatives $F'$ at each node, each element is mapped from reference space $\boldsymbol{\xi}$ to physical space $\boldsymbol{x}$ via a diffeomorphic map $\boldsymbol{x} = \boldsymbol{X}(\boldsymbol{\xi})$. The inverse of this map can then be used to map the element into the reference space defined to be $[0,1]^d$, where standard stencils can be used to calculate a derivative in $\boldsymbol{\xi}$ reference space. For example in 2-dimensions the gradient of a function $u$ can be calculated as
\begin{equation}
    \begin{pmatrix}
        \frac{\partial u}{\partial x} \\
        \frac{\partial u}{\partial y}
    \end{pmatrix}
    = 
    \underbrace{\begin{pmatrix}
        \frac{\partial \xi_1}{\partial x} &  \frac{\partial \xi_2}{\partial x}\\
        \frac{\partial \xi_1}{\partial y} &  \frac{\partial \xi_2}{\partial y}
    \end{pmatrix}}_{=(\nabla \boldsymbol{X})^{-1}}
    \begin{pmatrix}
        \frac{\partial u}{\partial \xi_1} \\
        \frac{\partial u}{\partial \xi_2}
    \end{pmatrix}
\end{equation}
As the derivatives in $\boldsymbol{\xi}$ reference space are taken over a tensor product domain, derivatives in each direction $\xi_i$ can be isolated to a line rather than having to use all nodes of the element. This means that only 1D stencils are needed in each direction $\xi_i$ to form the gradient in higher dimensions. Thus each node only uses the value of the flux $F$ at $d \cdot p$ other nodes to form the derivative $F'$. This overall construction is shown in Fig. \ref{fig:stencils2d}.

\subsubsection{Upwind construction}
At the boundary separating two or more elements an upwind stencil must be picked analogous to the one dimensional case. We describe a method for determining the upwind elements used to form the stencil and element boundary nodes which is based on the Petrov-Galerkin interpretation of the FUSE method outlined in Sect. \ref{sect:petrov}.

Recalling the Petrov-Galerkin framework, the model problem Eq. \ref{eq:conv_2d} can be written as
\begin{equation} \label{eq:petrov2d}
     \sum_j \int_\Omega \psi_i \phi_j dx \cdot \frac{\partial}{\partial t}u_j \ + \sum_d \sum_l \int_\Omega \psi_i \frac{\partial \phi_l}{\partial x^d} dx \cdot F_d(u_l) = 0
\end{equation}
where $F_d$ denotes the the $d$-th component of the flux $F$. As in 1D we pick the test functions $\psi_i$ to be equal to the nodal basis functions $\phi_i$ but with their support restricted only to upwind elements.

\begin{figure}
    \centering
    \includegraphics[scale=0.38]{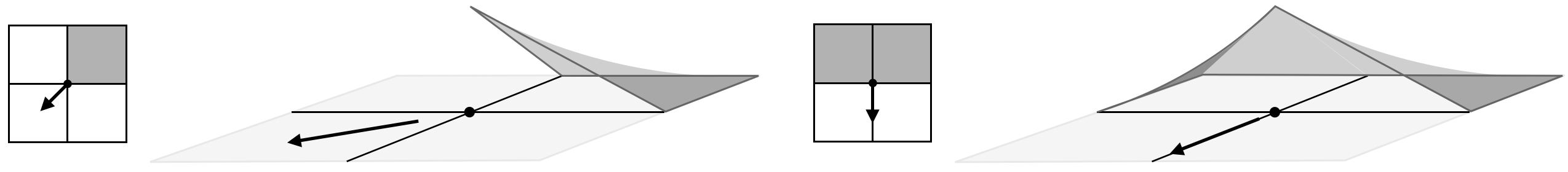}
    \caption{Schematic of upwind test functions in two dimensions. On the left the velocity at the centre node is pointing to the bottom left. In this case the top right element is upwind from the velocity and the test function at the node is simply the basis function at the centre node with its support limited to only that element. On the right the velocity is pointing straight down, so both elements on the top are upwind. In this case the test function at the node has support on both those elements.}
    \label{fig:upwind2d}
\end{figure}

To determine whether an element is upwind to a node $s_i$ on its boundary, the velocity $a\big(u(s_i)\big) = F'\big(u(s_i)\big)$ is first calculated at the node. Elements are upwind to the velocity at the node $s_i$ if the velocity can be traced backwards from the node into the element. Unlike in 1D this implies that more than one element may be upwind to a node, which may happen in the case that the velocity lines up with an element boundary. This procedure of determining upwind elements is shown in Fig. \ref{fig:upwind2d}.

\begin{figure}
    \centering
    \includegraphics[scale=0.35]{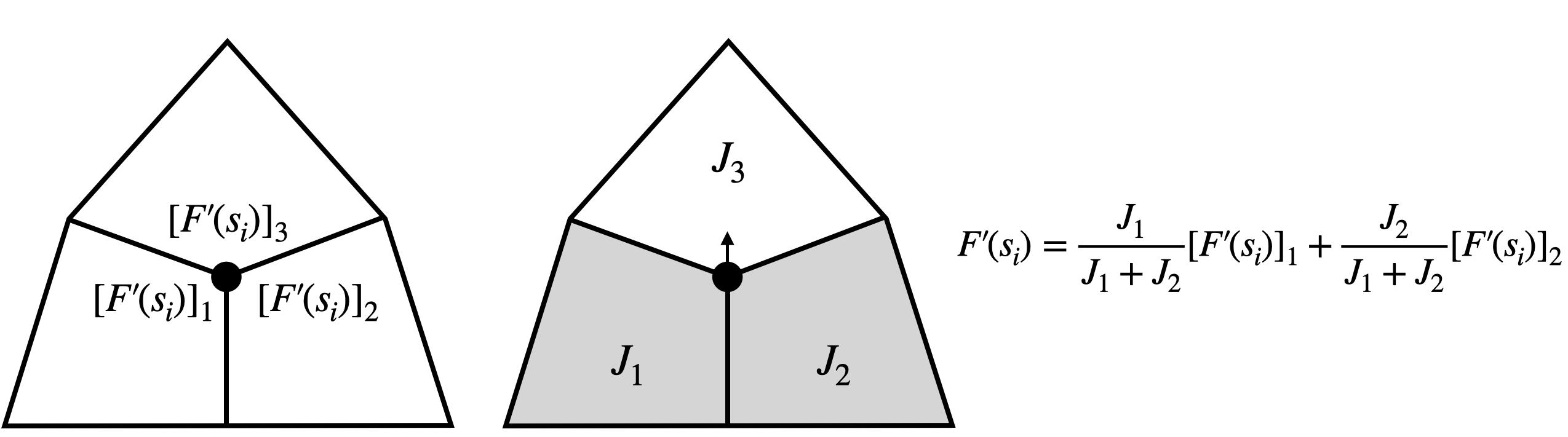}
    \caption{Schematic of upwind differencing in 2D. Centre node $s_i$ borders three elements and thus has a value for the derivative $F'$ from each element, denoted $[F'(s_i)]_j$. For a velocity then pointing straight up, the two elements on the bottom are both upwind so the value of the derivative at $s_i$ is given as a weighted average of the derivative values from these two elements.}
    \label{fig:upwind_stencils2d}
\end{figure}

To equate this with the spectral differencing formulation of the FUSE method, Eq. \ref{eq:petrov2d} is evaluated via nodal integration at the solution nodes. This allows for the following simplications as before using that: 1) $\phi_j(s_l) = \delta_{jl}$, 2) $\psi_i = \phi_i$ on the support of $\psi_i$ to get
\begin{equation}
      w_i \bigg( \frac{\partial}{\partial t}u_i \ + \sum_d \sum_l \mathbbm{1}_{\text{supp}(\psi_i)} \frac{d \phi_k^f}{dx^d} (s_i) \cdot F(u_l) \bigg) = 0
\end{equation}
similar to the 1D case. Denoting $K^{U}$ as the set of all upwind elements at the node $s_i$, this can then be written equivalently as
\begin{align}
    \frac{\partial}{\partial t}u_i \ &+ \sum_{K \in K^{U}} J_k \sum_{d,l} \frac{d \phi_k^f}{dx^d} (s_i) \cdot F(u_l) = 0 \\
    J_k &= \frac{1}{J} \int_{K} 1 ~dx, ~J = \sum_{K' \in K^{U}} \int_{K'} 1 ~dx
\end{align}
That is the derivative $F'$ at a node $s_i$ is given by an average of the spectral derivatives from all of the upwind elements weighted by their volume. A schematic of this is shown in Fig. \ref{fig:upwind_stencils2d}. 

\subsection{Relation to other methods}
\begin{figure}
    \centering
    \includegraphics[scale=0.35]{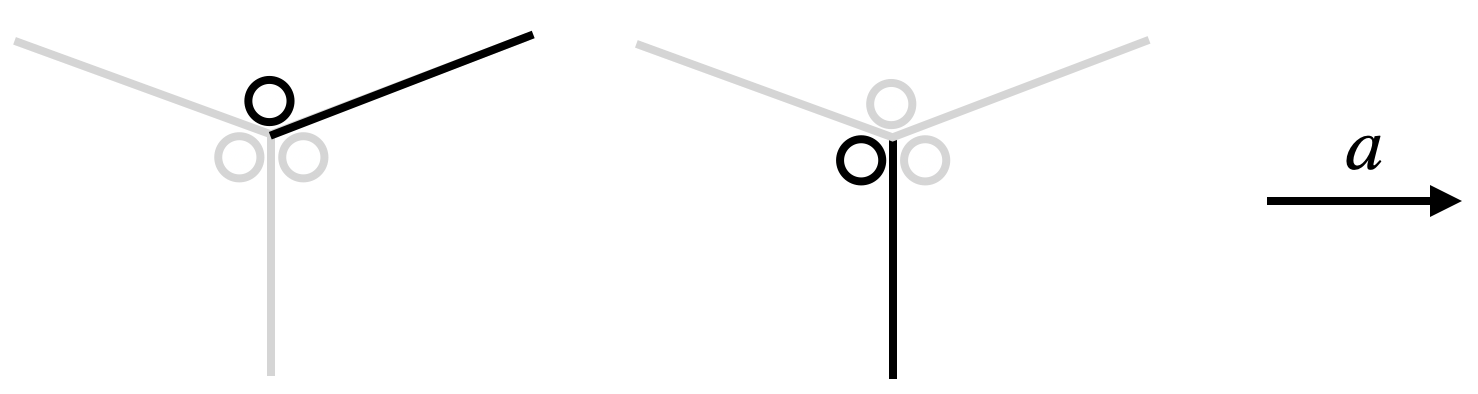}
    \caption{Example of SD method in 2D for constant-coefficient advection, where the velocity is pointing to the right. On the left for the highlighted edge, an upwind flux would take the value at the top highlighted node. On the right however for the highlighted edge an upwind flux would instead choose the value from the left highlighted node. As a result both these nodes are used in the SD method and can't be ignored to establish an equivalence with FUSE.}
    \label{fig:sd2d}
\end{figure}

Whilst in 2D the FUSE method remains equivalent to a nodally integrated Petrov-Galerkin method, unlike in 1D it can no longer be seen equivalent to a Spectral Differences even in the case of constant-coefficient advection. To see why this is the case we can consider a simple example in 2D with constant velocity pointing to the right; a schematic for this is shown in Fig. \ref{fig:sd2d}. The issue here is that at a corner separating multiple elements, as in SD and other related methods where a Riemann solver is applied on each boundary separating two elements, there is in general no way to pick a unique value for the flux at the corner. As a result the same construction in 1D where all repeated flux nodes except for one node at a corner separating more than two elements are ignored cannot be performed in higher dimensions.

\subsection{2nd derivative operators}
For second derivative operators in higher dimensions we again adapt the LDG method as in 1D. We note again that this is just one possible way of discretisation other methods such as Finite Elements possible given our choice of grid whereby nodes separating neighbouring elements are not duplicated.

We consider the same model problem in split form as in 1D
\begin{align}
    v &= \nabla u \\ 
    \frac{\partial u}{\partial t} - &\nabla \cdot v = 0 
\end{align}
To discretise the first equation a constant random velocity $a$ is chosen which is used to construct the upwind first order gradient operator. For the second equation then the divergence operator is discretised using the opposite velocity $-a$ such that it is the downwind analogue of the gradient operator. A random choice of velocity is used such that the probability is the velocity being along the direction of any line on the mesh is statistically zero. This results in a second derivative operator which has an upwind/downwind structure at each element boundary. For DG type discretisations it was shown in \cite{cockburn2007analysis} that this results in a stable second derivative operator and we observe this also for the FUSE method in numerical examples. 

\section{Numerical examples}
\subsection{Overview}
We apply the method to a selection of problems in one and two dimensions. Examples in three dimension are not included here for simplicity although the method can be easily extended to higher dimensions in space.

\subsection{Examples in 1D}
We first consider 1D examples demonstrating the two operators constructed in Sects. \ref{sect:conv1d} and \ref{sect:poi1d}. We are interested here both in the accuracy and spectral radii of the FUSE operators, as this gives a good estimate for the CFL number. We also verify that in all cases that the resulting operators are indeed stable.

In these examples for a given polynomial degree $p$, we look at the error in the discrete relative $L^\infty$ norm given by $\frac{\| u-u_0\|_\infty}{\| u_0 \|_\infty}$, where $u_0$ denotes the analytical solution, at the solution nodes and the spectral radii of the degree $p$ FUSE operator. To get a sense of the performance of the method we compare these with the corresponding degree $p-1$ DG operator for the considered problem on the same mesh. This choice is made as on a given mesh, the degree $p$ FUSE operator and degree $p-1$ DG operator have the same number of degrees of freedom (DOFs) due to the repeated boundary nodes in DG, which are not present in FUSE.

\subsubsection{Advection equation}
We consider the case of the advection equation in 1D
\begin{equation}
    \frac{\partial u}{\partial t} + \frac{\partial u}{\partial x} = 0
\end{equation}
on the domain $[0,1]$ with periodic boundary conditions. A uniform mesh is used to discretise the domain. To discretise the time derivative an RK4 integrator with a timestep of $10^{-4}$ is used. The final solution at $T=1$ is then checked against the initial solution which is set to a Gaussian pulse of $u_0(x) = \exp\big[-100 (x-0.5)^2 \big]$, where the solution is expected to be equal to the initial condition.

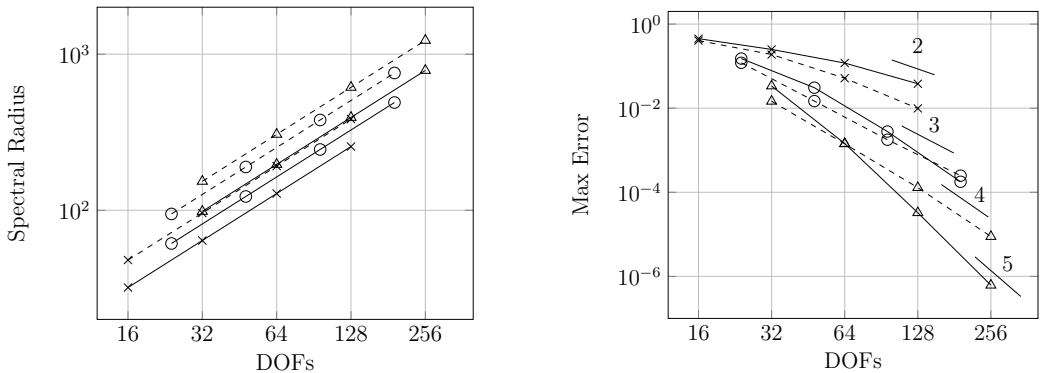
\begin{figure}[h]
    \centering
    \begin{minipage}{0.38\textwidth}
    \resizebox{\linewidth}{!}{%
    \begin{tikzpicture}[font=\large]
        \begin{loglogaxis}[
            xlabel = {DOFs},
            ylabel = {Spectral Radius},
            ylabel style={ yshift=2ex },
            ymin=2e1, ymax=2e3,
            xmin=12, xmax= 400,
            xtick={16,32,64,128,256},
            xticklabels={16,32,64,128,256},
            ytick={1e4,1e3,1e2},
            yticklabels={$10^{4}$,$10^{3}$,$10^{2}$},
            grid = both,
            grid style = {line width=.1pt, draw=gray!15},
            major grid style = {line width=.2pt, draw=gray!50},
        ]
        \addplot[mark=x, color=black, mark size=3pt]
        table[ x=n2, y=p2fuse ]{conv1d_cfl.txt};

        \addplot[dashed, mark=x, color=black, mark size=3pt, mark options={solid} ]
        table[ x=n2, y=p1dg ]{conv1d_cfl.txt};
        
        \addplot[mark=o, color=black, mark size=3pt]
        table[ x=n3, y=p3fuse ]{conv1d_cfl.txt};

        \addplot[dashed, mark=o, color=black, mark size=3pt, mark options={solid}]
        table[ x=n3, y=p2dg ]{conv1d_cfl.txt};

        \addplot[mark=triangle, color=black, mark size=3pt]
        table[ x=n4, y=p4fuse ]{conv1d_cfl.txt};

        \addplot[dashed, mark=triangle, color=black, mark size=3pt, mark options={solid}]
        table[ x=n4, y=p3dg ]{conv1d_cfl.txt};
        
        \end{loglogaxis}
    \end{tikzpicture}%
    }
    \end{minipage}
    \hspace{10mm}
    \begin{minipage}{0.38\textwidth}
    \resizebox{\linewidth}{!}{%
    \begin{tikzpicture}[font=\large]
        \begin{loglogaxis}[
            xlabel = {DOFs},
            ylabel = {Max Error},
            ylabel style={ yshift=2ex },
            ymin=1e-7, ymax=2e0,
            xmin=12, xmax= 400,
            xtick={16,32,64,128,256},
            xticklabels={16,32,64,128,256},
            ytick={1e-0,1e-2,1e-4,1e-6},
            yticklabels={$10^{0}$,$10^{-2}$,$10^{-4}$,$10^{-6}$},
            grid = both,
            grid style = {line width=.1pt, draw=gray!15},
            major grid style = {line width=.2pt, draw=gray!50},
        ]
        \addplot[mark=x, color=black, mark size=3pt]
        table[ x=n2, y=p2fuse ]{conv1d_comp.txt};

        \addplot[dashed, mark=x, color=black, mark size=3pt, mark options={solid} ]
        table[ x=n2, y=p1dg ]{conv1d_comp.txt};
        
        \addplot[mark=o, color=black, mark size=3pt]
        table[ x=n3, y=p3fuse ]{conv1d_comp.txt};

        \addplot[dashed, mark=o, color=black, mark size=3pt, mark options={solid}]
        table[ x=n3, y=p2dg ]{conv1d_comp.txt};

        \addplot[mark=triangle, color=black, mark size=3pt]
        table[ x=n4, y=p4fuse ]{conv1d_comp.txt};

        \addplot[dashed, mark=triangle, color=black, mark size=3pt, mark options={solid}]
        table[ x=n4, y=p3dg ]{conv1d_comp.txt};

        \addplot[domain=100:150, samples=2] {1400*2^(-2*log2(x))};
        \node at (axis cs:128, 3e-1) {2};

        \addplot[domain=110:180, samples=2] {5000*2^(-3*log2(x))};
        \node at (axis cs:150, 4e-3) {3};

        \addplot[domain=160:250, samples=2] {100000*2^(-4*log2(x))};
        \node at (axis cs:230, 9e-5) {4};

        \addplot[domain=220:340, samples=2] {1500000*2^(-5*log2(x))};
        \node at (axis cs:300, 2e-6) {5};
        
        \end{loglogaxis}
    \end{tikzpicture}%
    }
    \end{minipage}
    \caption{Advection equation example in 1D. We compare degree $p$ FUSE (solid lines) with degree $p-1$ DG (dashed lines) on the same mesh with equal number of refinements and therefore DOFs. Crosses show results for $p=2$, circles for $p=3$, and triangles for $p=4$. On the left the spectral radius for the advection operator is shown, and on the right the max error.}
    \label{fig:conv1d}
\end{figure}

The spectral radii and errors in the $L^\infty$ norm given are shown for several degrees $p$ on the left and right respectively in Fig. \ref{fig:conv1d}. For all degrees $p$, we observe that the spectral radius of degree $p$ FUSE to be around two-thirds of that of degree $p-1$ DG. In terms of accuracy, we find that FUSE for $p \geq 3$ seems to converge with order $O(h^{p+1})$, an extra order when compared to DG for the same number of DOFs. 

\subsubsection{Poisson equation}
For this example we solve Poisson's equation
\begin{equation}
    -\Delta u = f
\end{equation}
on the domain $[0,1]$ with Dirichlet boundary conditions. A uniform mesh is used also for this case. The right hand side $f$ is chosen such that the solution $u(x) = \exp\big[ \sin(2 \pi x) \big] - 1,$. 
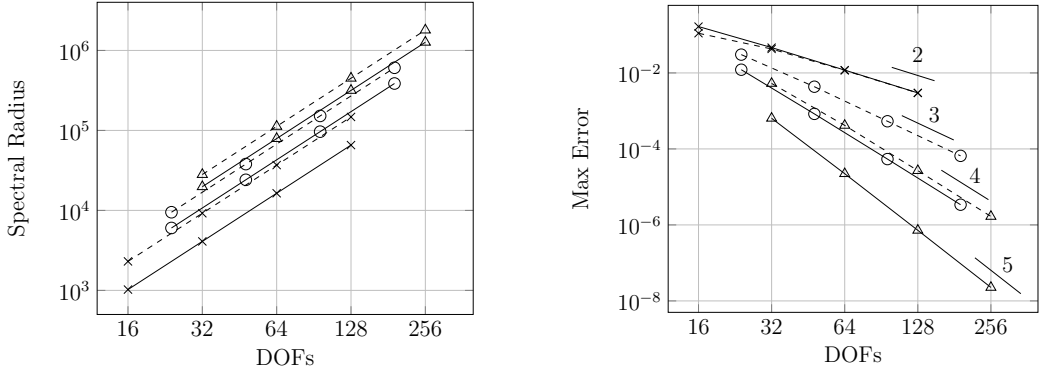
\begin{figure}[h]
    \centering
    \begin{minipage}{0.38\textwidth}
    \resizebox{\linewidth}{!}{%
    \begin{tikzpicture}[font=\large]
        \begin{loglogaxis}[
            xlabel = {DOFs},
            ylabel = {Spectral Radius},
            ylabel style={ yshift=2ex },
            ymin=5e2, ymax=4e6,
            xmin=12, xmax= 400,
            xtick={16,32,64,128,256},
            xticklabels={16,32,64,128,256},
            ytick={1e6,1e5,1e4,1e3},
            yticklabels={$10^{6}$,$10^{5}$,$10^{4}$,$10^{3}$},
            grid = both,
            grid style = {line width=.1pt, draw=gray!15},
            major grid style = {line width=.2pt, draw=gray!50},
        ]
        \addplot[mark=x, color=black, mark size=3pt]
        table[ x=n2, y=p2fuse ]{poi1d_cfl.txt};

        \addplot[dashed, mark=x, color=black, mark size=3pt, mark options={solid} ]
        table[ x=n2, y=p1dg ]{poi1d_cfl.txt};
        
        \addplot[mark=o, color=black, mark size=3pt]
        table[ x=n3, y=p3fuse ]{poi1d_cfl.txt};

        \addplot[dashed, mark=o, color=black, mark size=3pt, mark options={solid}]
        table[ x=n3, y=p2dg ]{poi1d_cfl.txt};

        \addplot[mark=triangle, color=black, mark size=3pt]
        table[ x=n4, y=p4fuse ]{poi1d_cfl.txt};

        \addplot[dashed, mark=triangle, color=black, mark size=3pt, mark options={solid}]
        table[ x=n4, y=p3dg ]{poi1d_cfl.txt};
        
        \end{loglogaxis}
    \end{tikzpicture}%
    }
    \end{minipage}
    \hspace{10mm}
    \begin{minipage}{0.38\textwidth}
    \resizebox{\linewidth}{!}{%
    \begin{tikzpicture}[font=\large]
        \begin{loglogaxis}[
            xlabel = {DOFs},
            ylabel = {Max Error},
            ylabel style={ yshift=2ex },
            ymin=5e-9, ymax=6e-1,
            xmin=12, xmax= 400,
            xtick={16,32,64,128,256},
            xticklabels={16,32,64,128,256},
            ytick={1e-2,1e-4,1e-6,1e-8},
            yticklabels={$10^{-2}$,$10^{-4}$,$10^{-6}$,$10^{-8}$},
            grid = both,
            grid style = {line width=.1pt, draw=gray!15},
            major grid style = {line width=.2pt, draw=gray!50},
        ]
        \addplot[mark=x, color=black, mark size=3pt]
        table[ x=n2, y=p2fuse ]{poi1d_comp.txt};

        \addplot[dashed, mark=x, color=black, mark size=3pt, mark options={solid} ]
        table[ x=n2, y=p1dg ]{poi1d_comp.txt};
        
        \addplot[mark=o, color=black, mark size=3pt]
        table[ x=n3, y=p3fuse ]{poi1d_comp.txt};

        \addplot[dashed, mark=o, color=black, mark size=3pt, mark options={solid}]
        table[ x=n3, y=p2dg ]{poi1d_comp.txt};

        \addplot[mark=triangle, color=black, mark size=3pt]
        table[ x=n4, y=p4fuse ]{poi1d_comp.txt};

        \addplot[dashed, mark=triangle, color=black, mark size=3pt, mark options={solid}]
        table[ x=n4, y=p3dg ]{poi1d_comp.txt};

        \addplot[domain=100:150, samples=2] {140*2^(-2*log2(x))};
        \node at (axis cs:128, 3e-2) {2};

        \addplot[domain=110:180, samples=2] {1000*2^(-3*log2(x))};
        \node at (axis cs:150, 8e-4) {3};

        \addplot[domain=160:250, samples=2] {18000*2^(-4*log2(x))};
        \node at (axis cs:220, 2e-5) {4};

        \addplot[domain=220:340, samples=2] {70000*2^(-5*log2(x))};
        \node at (axis cs:300, 9e-8) {5};
        
        \end{loglogaxis}
    \end{tikzpicture}%
    }
    \end{minipage}
    \caption{Poisson equation example in 1D. We compare degree $p$ FUSE (solid lines) with degree $p-1$ DG (dashed lines) on the same mesh with equal number of refinements and therefore DOFs. Crosses show results for $p=2$, circles for $p=3$, and triangles for $p=4$. On the left the spectral radius for the advection operator is shown, and on the right the max error.}
    \label{fig:poi1d}
\end{figure}

Fig. \ref{fig:poi1d} shows the spectral radius along with convergence plots in the $L^\infty$ norm for several degrees $p$. For the accuracy as with the previous example of the advection equation, for $p \geq 3$ degree $p$ FUSE is observed to converge at $p+1$ order, and an extra order when compared to degree $p-1$ DG with the same number of DOFs. However in terms of the spectral radius, while FUSE does for all orders shown here have a lower spectral radius, the amount by which it is lower decreases as $p$ increases. Significant gains in CFL number are therefore not expected in the case of diffusion dominated problems with FUSE over DG for high $p$.

\subsubsection{Euler equations}
As a final example we also consider the compressible Euler equations in 1D on the domain $\Omega = [0,1]$ with periodic boundary conditions. The Euler equations are given by
\begin{equation}
    \frac{\partial}{\partial t} 
    \underbrace{
    \begin{pmatrix}
        \rho \\
        \rho v \\
        \rho E
    \end{pmatrix}
    }_{u}
    +
    \frac{\partial}{\partial x}
    \underbrace{
    \begin{pmatrix}
        \rho u \\
        \rho v^2 + p \\
        v( \rho E + p )
    \end{pmatrix}}_{F(u)}
    =
    0
\end{equation}
where $\rho$ is the density, $v$ the velocity, $E$ the energy, and $p$ the pressure given by the equation of state $p = (\gamma - 1)(E-\frac{1}{2}\rho u^2)$. For this example, we set $\gamma = 1.4$, the initial condition to constant velocity and internal energy but with a smooth Gaussian for the density, and we run the simulation to a final time of $T=0.12$. To apply FUSE on this problem, we compute the derivatives of each component separately as before. For the upwinding, we rewrite the equations in quasilinear form
\begin{equation}
    \frac{\partial}{\partial t} u + A \cdot \frac{\partial}{\partial x} u = 0
\end{equation}
where $A = \frac{\partial F}{\partial u}$ is computed at the element boundaries. The tensor $A$ can be diagonalised as $A = \Lambda D \Lambda^{-1}$, allowing us to rewrite the equation in characteristic form
\begin{equation}
    \frac{\partial}{\partial t} \tilde u + D \frac{\partial}{\partial x} \tilde u = 0, ~\tilde{u} = \Lambda^{-1}u
\end{equation}
Each of the components of $\tilde{u}$ is known as a characteristic variable. As $D$ is diagonal, each characteristic variable is governed by a scalar conservation law of the form seen in Eq. \ref{eq:conserv1d}, and thus can be readily upwinded as with previous examples.

\begin{figure}[h]
    \centering
    \begin{minipage}{0.38\textwidth}
    \resizebox{\linewidth}{!}{%
    \begin{tikzpicture}[font=\large]
        \begin{loglogaxis}[
            xlabel = {DOFs per component},
            ylabel = {Max Error},
            ylabel style={ yshift=2ex },
            ymin=1e-8, ymax=2e0,
            xmin=30, xmax= 800,
            xtick={40,80,160,320,640},
            xticklabels={40,80,160,320,640},
            ytick={1e-0,1e-2,1e-4,1e-6,1e-8},
            yticklabels={$10^{0}$,$10^{-2}$,$10^{-4}$,$10^{-6}$,$10^{-8}$},
            grid = both,
            grid style = {line width=.1pt, draw=gray!15},
            major grid style = {line width=.2pt, draw=gray!50},
        ]
        \addplot[mark=x, color=black, mark size=3pt]
        table[ x=n2, y=p2fuse ]{euler1d_comp.txt};

        \addplot[dashed, mark=x, color=black, mark size=3pt, mark options={solid} ]
        table[ x=n2, y=p1dg ]{euler1d_comp.txt};
        
        \addplot[mark=o, color=black, mark size=3pt]
        table[ x=n3, y=p3fuse ]{euler1d_comp.txt};

        \addplot[dashed, mark=o, color=black, mark size=3pt, mark options={solid}]
        table[ x=n3, y=p2dg ]{euler1d_comp.txt};

        \addplot[mark=triangle, color=black, mark size=3pt]
        table[ x=n4, y=p4fuse ]{euler1d_comp.txt};

        \addplot[dashed, mark=triangle, color=black, mark size=3pt, mark options={solid}]
        table[ x=n4, y=p3dg ]{euler1d_comp.txt};

        \addplot[domain=330:600, samples=2] {3500*2^(-2*log2(x))};
        \node at (axis cs:450, 5e-2) {2};

        \addplot[domain=330:600, samples=2] {10000*2^(-3*log2(x))};
        \node at (axis cs:450, 3e-4) {3};

        \addplot[domain=330:600, samples=2] {180000*2^(-4*log2(x))};
        \node at (axis cs:450, 9e-6) {4};

        \addplot[domain=330:600, samples=2] {1500000*2^(-5*log2(x))};
        \node at (axis cs:450, 3.2e-8) {5};
        
        \end{loglogaxis}
    \end{tikzpicture}%
    }
    \end{minipage}
    \caption{Euler's equations example in 1D. We compare degree $p$ FUSE (solid lines) with degree $p-1$ DG (dashed lines) on the same mesh with equal number of refinements and therefore DOFs. Crosses show results for $p=2$, circles for $p=3$, and triangles for $p=4$.}
    \label{fig:euler1d}
\end{figure}
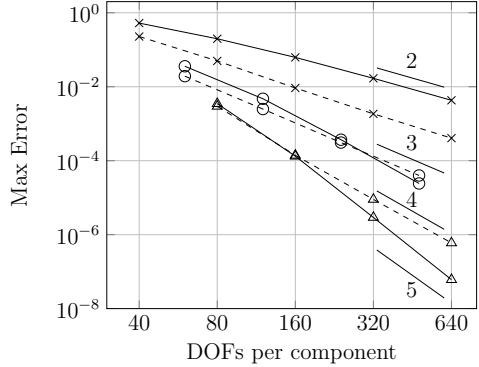

\subsection{Examples in 2D}
For the examples in 2D we consider a comparisons of the method on a structured versus unstructured mesh to obtain a measure of the method's performance on different geometries. Furthermore we demonstrate the method on a curvilinear mesh as which is well known to be necessary to obtain high-order accurate solutions on non-polygonal domains. Finally we consider an example of incompressible flow to demonstrate the method for CFD applications. For these examples we focus on the case $p=3$ as we feel it to be a good compromise between cost and accuracy whilst being high-order.

\subsubsection{Advection equation} \label{sect:advection2d}
We consider the advection equation in 2D
\begin{equation}
    \frac{\partial u}{\partial t} + \nabla \cdot (au) = 0
\end{equation}
on the domain $[0,1]^2$, with a divergence-free velocity chosen as $a(x,y) = (-y,x)$. This choice of velocity results in a velocity field that spins anticlockwise around the origin such that one travelling under this velocity field will return to their initial position after a time period of $2\pi$.

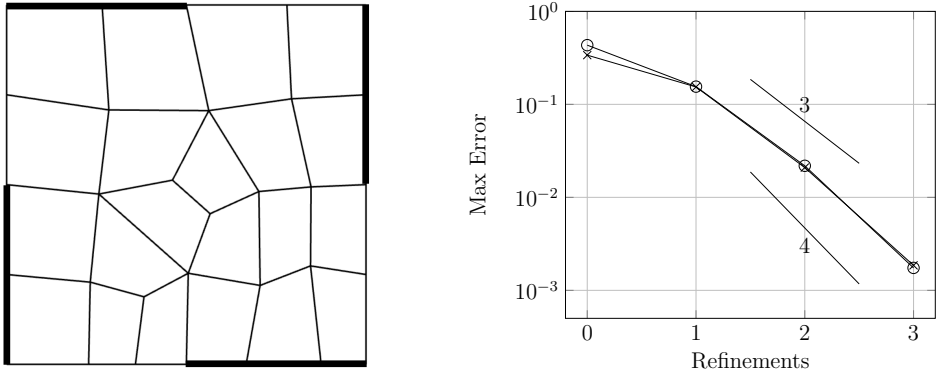
\begin{figure}[h]
    \centering
    \begin{minipage}{0.3\textwidth}
    \includegraphics[width=\textwidth]{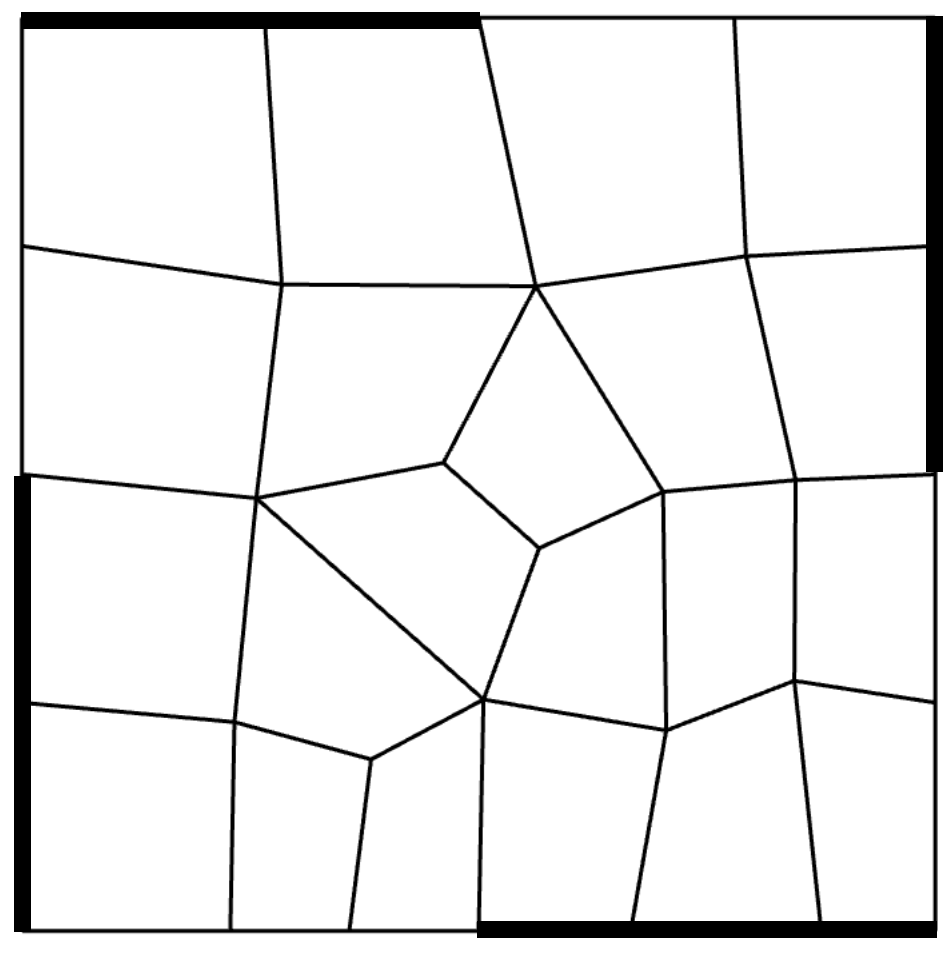}
    \end{minipage}
    \hspace{10mm}
    \begin{minipage}{0.38\textwidth}
    \resizebox{\linewidth}{!}{%
    \begin{tikzpicture}[font=\large]
        \begin{semilogyaxis}[
            xlabel = {Refinements},
            ylabel = {Max Error},
            ylabel style={ yshift=2ex },
            ymin=5e-4, ymax=1e-0,
            xmin=-0.2, xmax=3.2,
            xtick={0,1,2,3},
            xticklabels={0,1,2,3},
            ytick={1e-0, 1e-1,1e-2,1e-3,1e-4,1e-5},
            yticklabels={$10^{0}$,$10^{-1}$,$10^{-2}$,$10^{-3}$,$10^{-4}$,$10^{-5}$},
            grid = both,
            grid style = {line width=.1pt, draw=gray!15},
            major grid style = {line width=.2pt, draw=gray!50},
        ]
        \addplot[mark=x, color=black, mark size=3pt]
        table[ x=n, y=bad_p3 ]{conv2d_spin_p3.txt};
        \label{plot:conv2d_bad}
        
        \addplot[mark=o, color=black, mark size=3pt]
        table[ x=n, y=std_p3 ]{conv2d_spin_p3.txt};
        \label{plot:conv2d_std}

        \addplot[domain=1.5:2.5, samples=2] {4.2*2^(-3*x)};
        \node at (axis cs:2, 1e-1) {3};

        \addplot[domain=1.5:2.5, samples=2] {1.2*2^(-4*x)};
        \node at (axis cs:2, 3e-3) {4};
        
        \end{semilogyaxis}
    \end{tikzpicture}%
    }
    \end{minipage}
    \caption{Max error for 2D linear advection equation example on $[-1,1]^2$ for $p=3$. Unstructured mesh at zero refinements is shown on the left, where inflow boundaries for this problem have been bolded. \ref{plot:conv2d_bad} shows error on unstructured mesh, whilst \ref{plot:conv2d_std} shows error on structured mesh.}
    \label{fig:conv2d_spin_p3}
\end{figure}

The initial condition is defined as a Gaussian centred at $[-0.3,0]$ given by $u_0(x,y) = \exp\big[ -20 \big( (x+0.3)^2 + (y)^2 ) \big]$. The time derivative is as with the 1D advection example discretised using RK4 with a timestep of $10^{-4}$. The simulation is run until a final time of $2\pi$ at which point the solution is compared against the initial condition. Two choices of mesh are used in this example, a uniform mesh and an unstructured mesh at zero refinements shown on the left of Fig. \ref{fig:conv2d_spin_p3}, in which the bolded lines show the inflow boundary where zero inflow is imposed. The uniform mesh at zero refinements is chosen such that $h$ is as close as possible to the mean of lengths of all edges in the unstructured mesh whilst maintaining that $\frac{1}{h}$ is an integer.

The error in the $L^\infty$ norm given by $\max|u_0 - u|$ is shown on the right of Fig. \ref{fig:conv2d_spin_p3} for this example. Both display order of accuracy between $p$ and $p+1$ which does not seem to be affected by the presence of geometric defects on the unstructured mesh. We note however that while the discretisation on both meshes are stable that the of the operator on the unstructured mesh has a slightly larger spectral radius. This can likely be attributed to the fact that while the average edge lengths of both meshes are extremely comparable, the minimum edge length on the unstructured mesh is by design smaller than that on the structured one.

\subsubsection{Poisson equation}
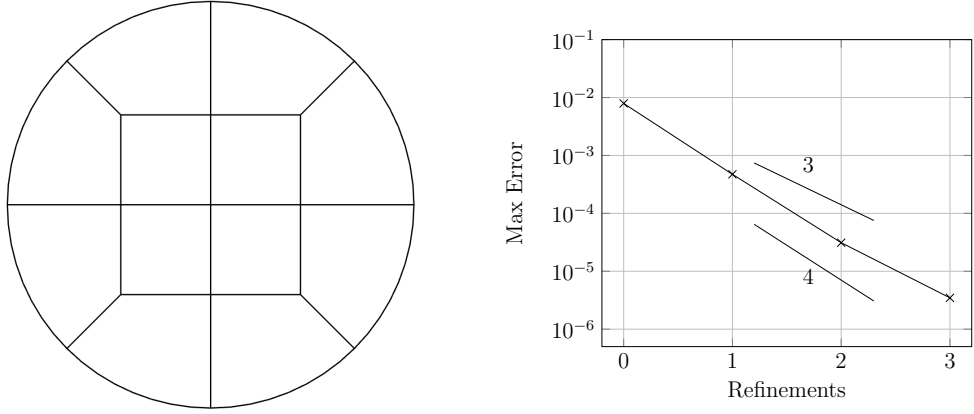
\begin{figure}[h]
    \centering
    \begin{minipage}{0.3\textwidth}
    \includegraphics[width=1.5\textwidth]{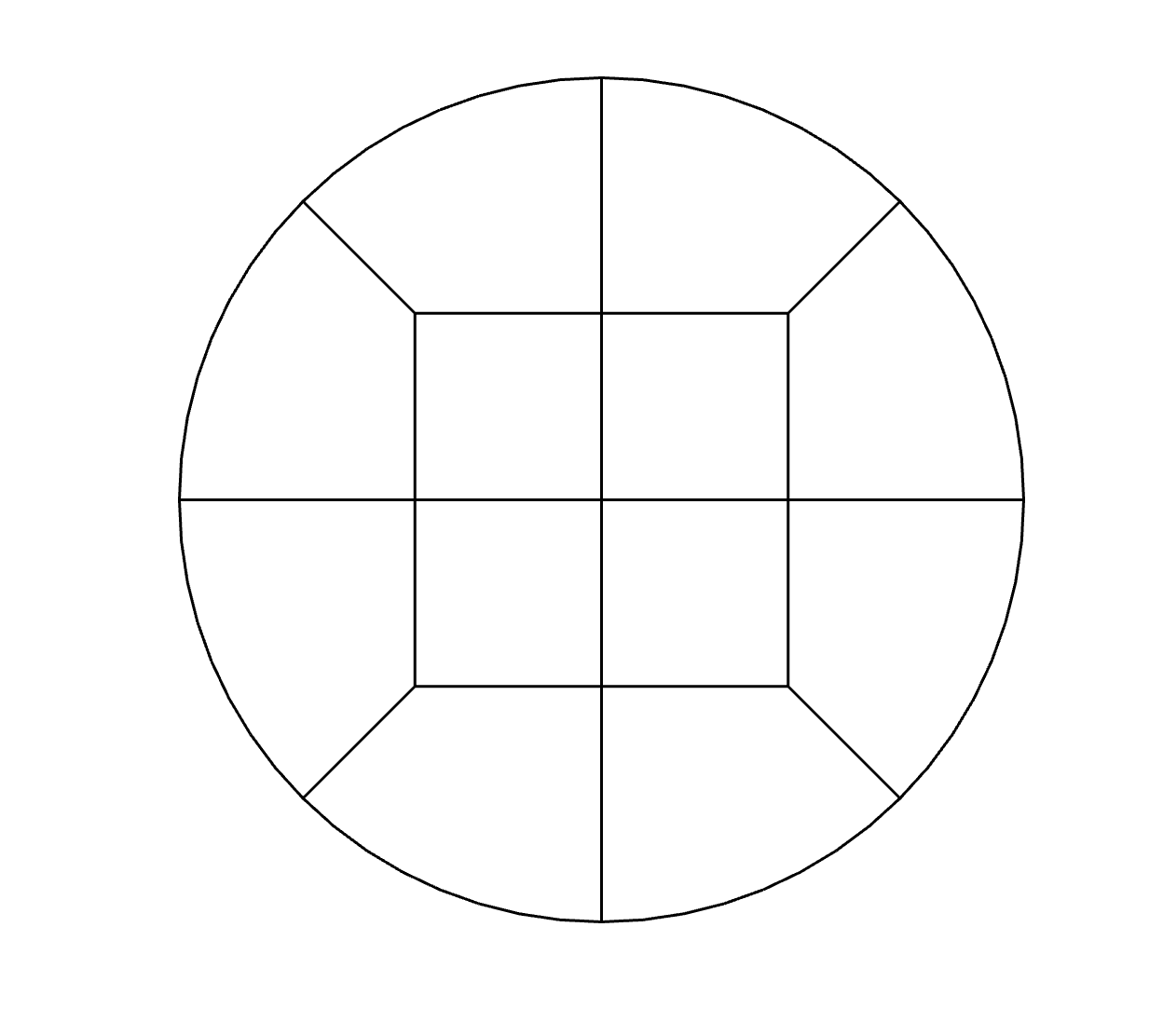}
    \end{minipage}
    \hspace{25mm}
    \begin{minipage}{0.38\textwidth}
    \resizebox{\linewidth}{!}{%
    \begin{tikzpicture}[font=\large]
        \begin{semilogyaxis}[
            xlabel = {Refinements},
            ylabel = {Max Error},
            ylabel style={ yshift=2ex },
            ymin=5e-7, ymax=1e-1,
            xmin=-0.2, xmax=3.2,
            xtick={0,1,2,3,4},
            xticklabels={0,1,2,3,4},
            ytick={1e-0, 1e-1,1e-2,1e-3,1e-4,1e-5,1e-6},
            yticklabels={$10^{0}$,$10^{-1}$,$10^{-2}$,$10^{-3}$,$10^{-4}$,$10^{-5}$,$10^{-6}$},
            grid = both,
            grid style = {line width=.1pt, draw=gray!15},
            major grid style = {line width=.2pt, draw=gray!50},
        ]
        \addplot[mark=x, color=black, mark size=3pt]
        table[ x=n, y=p3 ]{circpoi2d.txt};

        \addplot[domain=1.2:2.3, samples=2] {0.009*2^(-3*x)};
        \node at (axis cs:1.7, 7e-4) {3};

        \addplot[domain=1.2:2.3, samples=2] {0.0018*2^(-4*x)};
        \node at (axis cs:1.7, 8e-6) {4};
        
        \end{semilogyaxis}
    \end{tikzpicture}%
    }
    \end{minipage}
    \caption{Max error for 2D Poisson equation example on circle for $p=3$. Curved mesh at zero refinements is shown on the left. Error in $L^\infty$ norm as a function of number of refinements is shown on the right.}
    \label{fig:circpoi2d}
\end{figure}

We consider the Poisson equation
\begin{equation}
    -\Delta u = f
\end{equation}
on a circular domain centred at zero with radius $r=1$ as shown on the left in Fig. \ref{fig:circpoi2d}. The boundary of the domain is discretised using a cubic spline to match the degree $p$ used in the method. The right hand side $f$ and Dirichlet boundary conditions are chosen such that the solution $u(x,y) = \exp\big[ 1-x^2-y^2 \big]$.

The error in the $L^\infty$ norm is shown on the right of Fig. \ref{fig:circpoi2d}. We observe the order of accuracy to once again be between $p$ and $p+1$ consistent with the examples above. The accuracy of the method does not seem to be affected by the presence of curved boundaries allowing the method to be applied onto problems with complex geometries. We note however that with the presence of geometric defects in the mesh that whilst the eigenvalues of the discrete Laplace operator have positive real part the eigenvalues are in general complex.

\subsubsection{Incompressible flow}
We also consider the incompressible Navier-Stokes equations
\begin{align}
    \frac{\partial u}{\partial t} + (u \cdot \nabla) u - \nu &\Delta u + \frac{1}{\rho} \nabla p = 0 \\
    \nabla \cdot u &= 0
\end{align}
where $u, p$ denotes the velocity and pressure and constants $\rho, \nu$ are the density and kinematic viscosity of the system. The specific problem considered here is that of Taylor-Green vortex on the domain $[0,2\pi]^2$ with periodic boundary conditions. For this problem an analytical solution is known and given by:
\begin{align*}
    u_1 &= \sin(x) \cos(y) \exp(-2\nu t) \\
    u_2 &= \sin(x) \cos(y) \exp(-2\nu t) \\
    p &= \frac{\rho}{4} \big( \cos(2x) + \cos(2y) \big) \exp(-4\nu t)
\end{align*}
A structured mesh is used for this example where at the coarsest level each dimension is cut into four parts. The simulation is run to a final time of $T=0.1$ using a Crank-Nicolson time integrator with timestep $10^{-4}$, with problem constants set respectively to $\rho = 1, \nu = 1$.

The problem is discretised using the same LDG methodology as with the Poisson problem, whereby given a velocity field $v$ the discrete gradient operator $G_h$ is formed using the upwind construction described in Sect. \ref{sect:2d_1stderiv} and the discrete divergence operator $D_h$ formed as a downwind operator using the same velocity field $v$. The discrete Laplacian is as before discretised by composing the discrete divergence with the discrete gradient $L_h = D_h * G_h$. Introducing variables $u_h, p_h$ as the discrete solution approximations to the velocity and pressure the overall discretisation then reads as
\begin{align}
    \frac{\partial u}{\partial t} + (u_h \cdot G_h) u_h - \nu &D_h G_h u_h + \frac{1}{\rho} G_h p_h = 0 \\
    D_h \cdot u_h &= 0
\end{align}
The equation is solved for the discrete velocity and pressure at each timestep using Newton's method, with a starting guess of $(u_h,p_h) = (0,0)$ and iterated until the norm of the residual is no greater than $10^{-8}$.

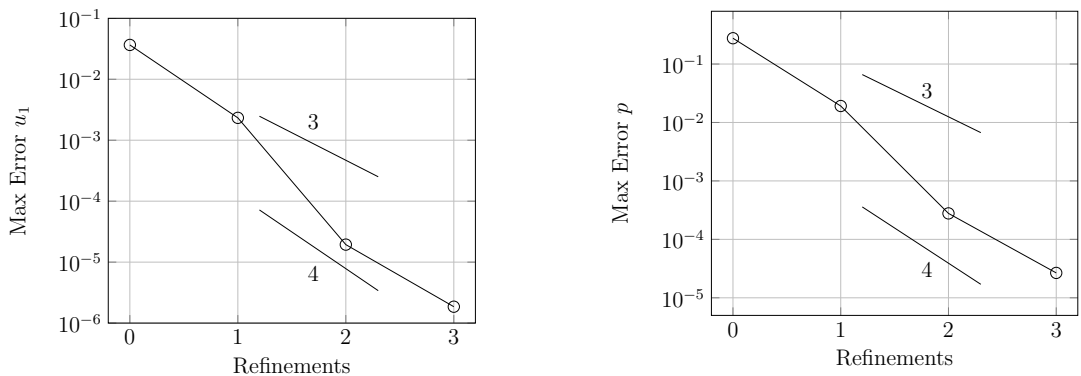
\begin{figure}[h]
    \centering
    \begin{minipage}{0.38\textwidth}
    \resizebox{\linewidth}{!}{%
    \begin{tikzpicture}[font=\large]
        \begin{semilogyaxis}[
            xlabel = {Refinements},
            ylabel = {Max Error $u_1$},
            ylabel style={ yshift=2ex },
            ymin=1e-6, ymax=1e-1,
            xmin=-0.2, xmax=3.2,
            xtick={0,1,2,3,4},
            xticklabels={0,1,2,3,4},
            ytick={1e0, 1e-1,1e-2,1e-3,1e-4,1e-5,1e-6},
            yticklabels={$10^{0-1}$,$10^{-1}$,$10^{-2}$,$10^{-3}$,$10^{-4}$, $10^{-5}$, $10^{-6}$},
            grid = both,
            grid style = {line width=.1pt, draw=gray!15},
            major grid style = {line width=.2pt, draw=gray!50},
        ]
        \addplot[mark=o, color=black, mark size=3pt]
        table[ x=n, y=structured_x ]{taylor_green.txt};

        \addplot[domain=1.2:2.3, samples=2] {3e-2*2^(-3*x)};
        \node at (axis cs:1.7, 2e-3) {3};

        \addplot[domain=1.2:2.3, samples=2] {2e-3*2^(-4*x)};
        \node at (axis cs:1.7, 6.5e-6) {4};
        
        \end{semilogyaxis}
    \end{tikzpicture}%
    }
    \end{minipage}
    \hspace{15mm}
    \begin{minipage}{0.38\textwidth}
    \resizebox{\linewidth}{!}{%
    \begin{tikzpicture}[font=\large]
        \begin{semilogyaxis}[
            xlabel = {Refinements},
            ylabel = {Max Error $p$},
            ylabel style={ yshift=2ex },
            ymin=5e-6, ymax=8e-1,
            xmin=-0.2, xmax=3.2,
            xtick={0,1,2,3,4},
            xticklabels={0,1,2,3,4},
            ytick={1e-0, 1e-1,1e-2,1e-3,1e-4,1e-5,1e-6},
            yticklabels={$10^{0}$,$10^{-1}$,$10^{-2}$,$10^{-3}$,$10^{-4}$,$10^{-5}$,$10^{-6}$},
            grid = both,
            grid style = {line width=.1pt, draw=gray!15},
            major grid style = {line width=.2pt, draw=gray!50},
        ]
        \addplot[mark=o, color=black, mark size=3pt]
        table[ x=n, y=structured_p ]{taylor_green.txt};

        \addplot[domain=1.2:2.3, samples=2] {0.8*2^(-3*x)};
        \node at (axis cs:1.8, 3.5e-2) {3};

        \addplot[domain=1.2:2.3, samples=2] {0.01*2^(-4*x)};
        \node at (axis cs:1.8, 3e-5) {4};
        
        \end{semilogyaxis}
    \end{tikzpicture}%
    }
    \end{minipage}
    \caption{Taylor-Green vortex example. On the left the max error in the first component of the velocity is shown, on the right the pressure.}
    \label{fig:taylorgreen}
\end{figure}

The error in the first component of the velocity and in the pressure are shown in Fig. \ref{fig:taylorgreen}. To obtain these plots the mesh was refined 3 times uniformly and polynomial degree $p=3$ was used. We observe that the $L^\infty$ error in both the velocity and pressure appears to converge with order of between $p$ and $p+1$ as with the previous examples. However while the results we obtain appear to be smooth and converge to the analytical solution for this example, it is unclear whether the method in general satisfies the the well known inf-sup stability conditions for incompressible flow. We plan to explore this in more detail in a future publication and do not consider it further here.

\section{Conclusion}
We have introduced a stabilised face-upwinded spectral element (FUSE) method for first and second order partial differential equations. The method is high order accurate and is suitable for use on general unstructured quadrilateral meshes. Nodes in this method are not duplicated across element boundaries akin to the Finite Element method resulting in fewer degrees of freedom. In particular it is stabilised via the specific choice of the Gauss-Legendre quadarature points plus endpoints node distribution in each element as well as picking upwinded stencils on element boundaries.

We plan in a future work to study in more detail the behaviour of the method for incompressible flows and in particular on its stability properties for this type of problem. Furthermore we intend to explore in more depth the performance of the method on more complex problems and its comparison with other popular stabilised methods in practice. Finally we also remain interested in extending the method to meshes of element shapes other than quadrilaterals, most importantly in constructing a stable extension for simplex based meshes as is used in many commercial softwares.

\section*{Acknowledgments}

This work was supported in part by the Director, Office of Science, Office of Advanced Scientific Computing Research, U.S. Department of Energy under Contract No. DE-AC02-05CH11231.

\newpage
\bibliographystyle{plain}
\bibliography{references}

\end{document}